\documentclass{article}[12pt]

\usepackage{latexsym}
\usepackage{amsmath}
\usepackage{amssymb}

\begin{document}

\title{Solvable systems of two coupled first-order ODEs with homogeneous
cubic polynomial right-hand sides}

\author{Francesco Calogero$^{a,b}$\thanks{e-mail: francesco.calogero@roma1.infn.it}
\thanks{e-mail: francesco.calogero@uniroma1.it}
 , Farrin Payandeh$^c$\thanks{e-mail: farrinpayandeh@yahoo.com}
 \thanks{e-mail: f$\_$payandeh@pnu.ac.ir}}

\maketitle   \centerline{\it $^{a}$Physics Department, University of
Rome "La Sapienza", Rome, Italy}

\maketitle   \centerline{\it $^{b}$INFN, Sezione di Roma 1}

\maketitle

\maketitle   \centerline{\it $^{c}$Department of Physics, Payame
Noor University, PO BOX 19395-3697 Tehran, Iran}

\maketitle

\begin{abstract}

The solution $x_{n}\left( t\right) ,$ $n=1,2,$ of the \textit{initial-values}
problem is reported of the \textit{autonomous} system of $2$ coupled
first-order ODEs with \textit{homogeneous cubic polynomial} right-hand sides,%
\begin{eqnarray}
\dot{x}_{n}=c_{n1}\left( x_{1}\right) ^{3}+c_{n2}\left( x_{1}\right)
^{2}x_{2}+c_{n3}x_{1}\left( x_{2}\right) ^{2}+c_{n4}\left( x_{2}\right)
^{3}~,~~~n=1,2~,   \nonumber
\end{eqnarray}%
when the $8$ (time-independent) coefficients $c_{n\ell }$ are appropriately
defined in terms of $7$ \textit{arbitrary} parameters, which then also
identify the solution of this model. The inversion of these relations is
also investigated, namely how to obtain, in terms of the $8$ coefficients $%
c_{n\ell },$ the $7$ parameters characterizing the solution of this model;
and $2$ \textit{constraints} are \textit{explicitly} identified which, if
satisfied by the $8$ parameters $c_{n\ell },$ guarantee the \textit{%
solvability by algebraic operations} of this dynamical system. Also
identified is a related, \textit{appropriately modified}, class of (generally
\textit{complex}) systems, reading%
\begin{eqnarray}
\dot{{\tilde{x}}_{n}}=\mathbf{i}\omega \tilde{x}_{n}+c_{n1}\left(
\tilde{x}_{1}\right) ^{3}+c_{n2}\left( \tilde{x}_{1}\right) ^{2}\tilde{x}%
_{2}+c_{n3}\tilde{x}_{1}\left( \tilde{x}_{2}\right) ^{2}+c_{n4}\left( \tilde{%
x}_{2}\right) ^{3}~,~~~n=1,2~, \nonumber
\end{eqnarray}%
with $\mathbf{i}\omega $ an \textit{arbitrary imaginary} parameter, which
feature the remarkable property to be \textit{isochronous}, namely their
\textit{generic} solutions are---as functions of \textit{real time}---%
\textit{completely periodic} with a period which is, for each of these
models, a \textit{fixed} \textit{integer multiple} of the basic period $%
\tilde{T}=2\pi /\left\vert \omega \right\vert $.

\end{abstract}

\section{Introduction and presentation of the main results}

The system characterized by the following $2$ \textit{nonlinearly-coupled}
Ordinary Differential Equations (ODEs) with \textit{homogeneous cubic
polynomial} right-hand sides,
\begin{equation}
\dot{x}_{n}=c_{n1}\left( x_{1}\right) ^{3}+c_{n2}\left( x_{1}\right)
^{2}x_{2}+c_{n3}x_{1}\left( x_{2}\right) ^{2}+c_{n4}\left( x_{2}\right)
^{3}~,~~~n=1,2~,  \label{1}
\end{equation}%
is a prototypical example of \textit{autonomous} dynamical systems.

\textbf{Notation 1-1}. The $2$ (possibly \textit{complex}) numbers $%
x_{n}\equiv x_{n}\left( t\right) ,$ $n=1,2,$ are the dependent variables; $t$
is the independent variable ("time"; but the treatment remains valid when $t$
is considered a \textit{complex} variable); superimposed dots denote $t$%
-differentiations; the $8$ $t$-independent (possibly \textit{complex})
coefficients $c_{n\ell }$ ($n=1,2;$ $\ell =1,2,3,4$) can be \textit{a priori
arbitrarily assigned}; but a key result of this paper is to identify $2$
\textit{algebraic constraints} (see below \textbf{Subsection 3.2} and
\textbf{Appendix C}) which---if satisfied by the $8$ coefficients $c_{n\ell
} $---allows the solution of the \textit{initial-values} problem of the
system (\ref{1}) \textit{by algebraic operations }(see below \textbf{%
Proposition 1-1}). $\blacksquare $

\textbf{Remark 1-1}. It is plain that the system (\ref{1}) is \textit{%
invariant} under the following simultaneous rescaling of the $8$
coefficients $c_{n\ell }$ and of the independent variable $t$:%
\begin{equation}
x_{n}\left( t\right) \Rightarrow \eta x_{n}\left( \eta ^{2}t\right)
~,~~~n=1,2~,  \label{1Inveta}
\end{equation}%
with $\eta $ an \textit{arbitrary nonvanishing parameter}.

It is also plain (see (\ref{1})) that, by \textit{rescaling} each of the $2$
dependent variables $x_{1}\left( t\right) $ and $x_{2}\left( t\right) $ by
appropriate \textit{constant} parameters, $2$ of the $8$ parameters $%
c_{n\ell }$ may generally be replaced by \textit{arbitrary} constants (for
instance by \textit{unity}; of course unless that coefficient \textit{%
vanishes} to begin with); while the other coefficients are of course also
appropriately \textit{rescaled}. But of course these features are of no help
to satisfy the $2$ \textit{constraints} on the $8$ coefficients $c_{n\ell }$%
---shown below to be sufficient for the \textit{algebraic solvability} of
the model (\ref{1})---which are \textit{invariant} under such \textit{%
rescalings} (see below). In the following we generally assume that the
parameters $c_{n\ell }$ have \textit{generic} values (except for satisfying
the \textit{constraints} mentioned just above); the only instance violating
this rule is the special case (see below \textbf{Subsection 3.5}) with $%
c_{14}=c_{21}=0$ which is particularly relevant in \textit{applicative}
contexts and deserve therefore a separate treatment. $\blacksquare $

\textbf{Remark 1-2}. It is plain that the system (\ref{1}) is \textit{%
invariant} under the following exchange of variables and parameters:%
\begin{equation}
x_{1}\left( t\right) \Leftrightarrow x_{2}\left( t\right)
~;~c_{11}\Leftrightarrow c_{24}~,~c_{12}\Leftrightarrow
c_{23}~,~c_{13}\Leftrightarrow c_{22}~,~c_{14}\Leftrightarrow c_{21}~.
\label{1Change}
\end{equation}%
This property shall be used occasionally below to decrease the number of
displayed equations. $\blacksquare $

\textbf{Remark 1-3}. The preceding \textbf{Remark 1-2 }might motivate the
\textit{merely notational} replacement---possibly preferable in some \textit{%
applicative} contexts---of the system (\ref{1}) with%
\begin{eqnarray}
\dot{\xi} &=&\alpha _{1}\xi ^{3}+\alpha _{2}\xi ^{2}\eta +\alpha _{3}\xi
\eta ^{2}+\alpha _{4}\eta ^{3}~,  \nonumber \\
\dot{\eta} &=&\beta _{1}\eta ^{3}+\beta _{2}\eta ^{2}\xi +\beta _{3}\xi \eta
^{2}+\beta _{4}\xi ^{3}~,
\end{eqnarray}%
corresponding to the change of notation $x_{1}\left( t\right) \Rightarrow
\xi \left( t\right) ,$ $x_{2}\left( t\right) \Rightarrow \eta \left(
t\right) ,$ $c_{1\ell }\Rightarrow \alpha _{\ell },$ $c_{2\ell }\Rightarrow
\beta _{5-\ell },$ $\ell =1,2,3,4$ and implying that the invariance property
(\ref{1Change}) take the neater form $\xi \left( t\right) \Leftrightarrow
\eta \left( t\right) ,$ $\alpha _{\ell }\Leftrightarrow \beta _{\ell },$ $%
\ell =1,2,3,4$. But throughout this paper we stick with the notation implied
by the more common notation of the system (\ref{1}). $\blacksquare $

The system (\ref{1}) has been investigated over time in an enormous number
of purely \textit{mathematical}, or mainly \textit{applicative}, contexts.
The first approach generally focussed on \textit{qualitative} features of
its solutions $x_{n}\left( t\right) $ as functions of \textit{real} $t$
("time"): mainly on the identification of its equilibria (if any), and on
the behaviors close to them and asymptotically at large times ($t\rightarrow
\infty $); and, as functions of \textit{complex} $t,$ on the \textit{%
analyticity} properties of its \textit{general solution}. The second
approach was motivated by \textit{applications} in various contexts (mainly
the time evolution of interacting "populations" or of other quantifiable
entities, such as concentrations of chemicals or financial entries, you name
it); and it was often pursued via numerical computations. This literature is
too large to allow any attempt to provide a list of references which would
do justice to the multitude of relevant papers. Here we limit ourselves to
identify only the relatively recent (open access) paper \cite{GRN2016}
containing several references, and the $2$ quite recent papers \cite{CP2019}
and \cite{CCL2020} because they have motivated this research and because
from them additional relevant references can be traced. But we also like to
mention the fundamental paper \cite{RG1960} by Ren\'{e} Garnier, who 60
years ago investigated this type of systems---in fact, starting from the
more general class of system analogous to (\ref{1}) but featuring \textit{a
priori arbitrary} polynomials on their right-hand sides and identifying
subcases of this general class (including subcases of eq. (\ref{1})) such
that their \textit{general} solutions---as functions of the independent
variable $t$ considered as a \textit{complex} variable---are \textit{%
holomorphic} ("uniform" in Garnier's language). The present paper may be
seen as a generalization of (some aspects) of Garnier's work, as it allows
the identification of a much larger subclass of the system (\ref{1}), those
which are solvable by \textit{algebraic operations}; hence such that their
\textit{general} \textit{solutions} (\textit{rather} \textit{explicitly}
identified below) only feature---as functions of the independent variable $t$
considered as a \textit{complex} variable---a \textit{finite} number of
branch points (i. e., singularities of type $\left( t-t_{s}\right) ^{r_{s}}$%
). This much larger subclass is clearly interesting from a purely
mathematical point of view, and even more so because---as detailed below,
see \textbf{Section 4}---the simple (\textit{complex}) extension of the
model (\ref{1}) reading
\begin{equation}
\dot {{\tilde{x}}_{n}}=\mathbf{i}\omega \tilde{x}_{n}+c_{n1}\left(
\tilde{x}_{1}\right) ^{3}+c_{n2}\left( \tilde{x}_{1}\right) ^{2}\tilde{x}%
_{2}+c_{n3}\tilde{x}_{1}\left( \tilde{x}_{2}\right) ^{2}+c_{n4}\left( \tilde{%
x}_{2}\right) ^{3}~,~~~n=1,2  \label{1Iso}
\end{equation}%
(with $\mathbf{i}\omega $ an \textit{arbitrary} additional, \textit{purely
imaginary}, parameter) features then---provided all the exponents $r_{s}$
are \textit{rational} numbers (a restriction that can be easily imposed, see
below)---the remarkable property of \textit{isochrony}, its \textit{generic}
solutions being then \textit{completely periodic} with a period which is an
\textit{integer multiple} of the basic period $\tilde{T}=2\pi /\left\vert
\omega \right\vert $ (see for instance, \cite{FC2008}). $\blacksquare $

Let us now state our main finding.

\textbf{Proposition 1-1}. Assume that the $8$ parameters $c_{n\ell }$ of the
dynamical system (\ref{1}) are given in terms of the $7$, \textit{a priori
arbitrary}, parameters $a_{1},$ $a_{2},~b_{1},$ $b_{2},$ $\gamma _{1},$ $%
\gamma _{2},$ $\gamma _{3}$ by the following $8$ formulas:
\begin{subequations}
\label{2}
\begin{equation}
c_{11}=\left[ \left( a_{1}\right) ^{3}b_{2}-a_{2}K_{1}\right] /c~,
\label{1c11}
\end{equation}%
\begin{equation}
c_{21}=\left[ -\left( a_{1}\right) ^{3}b_{1}+a_{1}K_{1}\right] /c~,
\label{1c21}
\end{equation}%
\begin{equation}
c_{12}=\left[ 3\left( a_{1}\right) ^{2}a_{2}b_{2}-a_{2}K_{2}\right] /c~,
\label{1c12}
\end{equation}%
\begin{equation}
c_{22}=\left[ -3\left( a_{1}\right) ^{2}a_{2}b_{1}+a_{1}K_{2}\right] /c~,
\label{1c22}
\end{equation}%
\begin{equation}
c_{13}=\left[ 3a_{1}\left( a_{2}\right) ^{2}b_{2}-a_{2}K_{3}\right] /c~,
\label{1c13}
\end{equation}%
\begin{equation}
c_{23}=\left[ -3a_{1}\left( a_{2}\right) ^{2}b_{1}+a_{1}K_{3}\right] /c~,
\label{1c23}
\end{equation}%
\begin{equation}
c_{14}=\left[ \left( a_{2}\right) ^{3}b_{2}-a_{2}K_{4}\right] /c~,
\label{1c14}
\end{equation}%
\begin{equation}
c_{24}=\left[ -\left( a_{2}\right) ^{3}b_{1}+a_{1}K_{4}\right] /c~,
\label{1c24}
\end{equation}%
where, above and hereafter,
\end{subequations}
\begin{subequations}
\label{1KK}
\begin{equation}
K_{1}\equiv \gamma _{1}a_{1}\left( b_{1}\right) ^{2}+\gamma _{2}\left(
a_{1}\right) ^{2}b_{1}+\gamma _{3}\left( a_{1}\right) ^{3}+\left(
b_{1}\right) ^{3}~,  \label{K1}
\end{equation}%
\begin{eqnarray}
K_{2} &\equiv &\gamma _{1}\left[ a_{2}\left( b_{1}\right)
^{2}+2a_{1}b_{1}b_{2}\right] +\gamma _{2}\left[ \left( a_{1}\right)
^{2}b_{2}+2a_{1}a_{2}b_{1}\right] +3\gamma _{3}\left( a_{1}\right) ^{2}a_{2}
\nonumber \\
&&+3\left( b_{1}\right) ^{2}b_{2}~,  \label{K2}
\end{eqnarray}%
\begin{eqnarray}
K_{3} &\equiv &\gamma _{1}\left[ a_{1}\left( b_{2}\right)
^{2}+2a_{2}b_{1}b_{2}\right] +\gamma _{2}\left[ \left( a_{2}\right)
^{2}b_{1}+2a_{1}a_{2}b_{2}\right] +3\gamma _{3}a_{1}\left( a_{2}\right) ^{2}
\nonumber \\
&&+3b_{1}\left( b_{2}\right) ^{2}~,  \label{K3}
\end{eqnarray}%
\begin{equation}
K_{4}\equiv \gamma _{1}a_{2}\left( b_{2}\right) ^{2}+\gamma _{2}\left(
a_{2}\right) ^{2}b_{2}+\gamma _{3}\left( a_{2}\right) ^{3}+\left(
b_{2}\right) ^{3}~,  \label{K4}
\end{equation}%
\begin{equation}
c\equiv a_{1}b_{2}-a_{2}b_{1}~.  \label{cc}
\end{equation}

Then the \textit{initial-values} problem for the system (\ref{1}) is \textit{%
solved} by the following formulas:
\end{subequations}
\begin{subequations}
\label{x1x2wt}
\begin{eqnarray}
x_{1}\left( t\right) &=&\left[ b_{2}y\left( t\right) -a_{2}w\left( t\right) %
\right] /c~,  \nonumber \\
x_{2}\left( t\right) &=&-\left[ b_{1}y\left( t\right) -a_{1}w\left( t\right) %
\right] /c~,  \label{x12t}
\end{eqnarray}%
with $c$ defined---above and hereafter---by eq. (\ref{cc}). Note that these
formulas are easily inverted, reading then
\begin{equation}
y\left( t\right) =a_{1}x_{1}\left( t\right) +a_{2}x_{2}\left( t\right)
,~~~w\left( t\right) =b_{1}x_{1}\left( t\right) +b_{2}x_{2}\left( t\right) ~,
\label{ywxt}
\end{equation}%
which clearly also imply
\begin{equation}
y\left( 0\right) =a_{1}x_{1}\left( 0\right) +a_{2}x_{2}\left( 0\right)
,~~~w\left( 0\right) =b_{1}x_{1}\left( 0\right) +b_{2}x_{2}\left( 0\right) ~.
\label{y0w0}
\end{equation}%
While the $2$ functions $y\left( t\right) $ and $w\left( t\right) $ are
given in terms of the above parameters and of their initial values $y\left(
0\right) $ and $w\left( 0\right) $ (themselves given in terms of the initial
data $x_{1}\left( 0\right) $ and $x_{2}\left( 0\right) $ by (\ref{y0w0})) by
the following \textit{explicit} formula for $y\left( t\right) $,
\end{subequations}
\begin{equation}
y\left( t\right) =y\left( 0\right) \left\{ 1-2\left[ y\left( 0\right) \right]
^{2}t\right\} ^{-1/2}~,  \label{Solyt}
\end{equation}%
and by the relation
\begin{subequations}
\label{Solwt}
\begin{equation}
w\left( t\right) =u\left( t\right) y\left( t\right) ~,  \label{wwwt}
\end{equation}%
where the function $u\left( t\right) $ is defined \textit{implicitly} by the
following formula:
\begin{equation}
\left[ \frac{u\left( t\right) -u_{1}}{u\left( 0\right) -u_{1}}\right]
^{-2\lambda _{1}}\left[ \frac{u\left( t\right) -u_{2}}{u\left( 0\right)
-u_{2}}\right] ^{-2\lambda _{2}}\left[ \frac{u\left( t\right) -u_{3}}{%
u\left( 0\right) -u_{3}}\right] ^{-2\lambda _{3}}=1-2\left[ y\left( 0\right) %
\right] ^{2}t~.  \label{1ut}
\end{equation}%
Here of course
\begin{equation}
u\left( 0\right) =w\left( 0\right) /y\left( 0\right) ~;  \label{1u0}
\end{equation}%
while the $6$ parameters $u_{j}$ and $\lambda _{j}$ ($j=1,2,3$) are
defined---in terms of (only!) the $3$ parameters $\gamma _{1},$ $\gamma
_{2}, $ $\gamma _{3}$---by the following \textit{cubic} equation,
\end{subequations}
\begin{subequations}
\label{alphauj}
\begin{equation}
u^{3}+\gamma _{1}u^{2}+\left( \gamma _{2}-1\right) u+\gamma _{3}=\left(
u-u_{1}\right) \left( u-u_{2}\right) \left( u-u_{3}\right) ~,
\end{equation}%
where $u$ is an arbitrary variable. The $3$ roots $u_{j}$ of this cubic
polynomial are of course related to the $3$ parameters $\gamma _{j}$ as
follows,
\begin{eqnarray}
\gamma _{1} &=&-\left( u_{1}+u_{2}+u_{3}\right) ~,  \nonumber \\
\gamma _{2} &=&u_{1}u_{2}+u_{2}u_{3}+u_{3}u_{1}+1~,  \nonumber \\
\gamma _{3} &=&-u_{1}u_{2}u_{3}~;
\end{eqnarray}%
and are explicitly given in terms of the $3$ parameters $\gamma _{j}$ by the
standard Cardano formulas. And they determine the $3$ parameters $\lambda
_{j}$ via the following $3$ relations,
\end{subequations}
\begin{subequations}
\label{landuj}
\begin{eqnarray}
\lambda _{1}+\lambda _{2}+\lambda _{3} &=&0~,  \nonumber \\
\lambda _{1}u_{1}+\lambda _{2}u_{2}+\lambda _{3}u_{3} &=&0~,  \nonumber \\
\lambda _{1}u_{2}u_{3}+\lambda _{2}u_{3}u_{1}+\lambda _{3}u_{1}u_{2} &=&1~,
\end{eqnarray}%
implying
\begin{equation}
\lambda _{j}=\left[ \prod\limits_{\ell =1,~\ell \neq j}^{3}\left(
u_{j}-u_{\ell }\right) \right] ^{-1}~,~~~j=1,2,3~.~~~\blacksquare
\label{landan}
\end{equation}

\textbf{Remark 1-4}. The relation (\ref{1ut}) implies that \textit{if} the $%
3 $ parameters $\lambda _{j},$ $j=1,2,3,$ are $3$ \textit{arbitrary real
rational} numbers,
\end{subequations}
\begin{equation}
\lambda _{j}=N_{j}/M_{j}~,~~~j=1,2,3
\end{equation}%
(with $N_{j}$ \textit{arbitrary integers} and $M_{j}$ \textit{arbitrary
positive integers}), then $u\left( t\right) $ is an \textit{algebraic}
function of $t,$ being then the root of the \textit{polynomial} equation
determining $u\left( t\right) $---namely the equation that obtains by first
multiplying eq. (\ref{1ut}) by the factor(s) $\left[ u\left( t\right) -u_{1}%
\right] ^{2\lambda _{j}}$ with $\lambda _{j}>0$ and then raising the
resulting equation to an appropriate positive integer power. Then the
solutions $x_{1}\left( t\right) $ and $x_{2}\left( t\right) $ are clearly as
well \textit{algebraic} functions of the time $t$. A remarkable consequence (%
\textit{isochrony}!) of this fact has already been mentioned above (see (\ref%
{1Iso})) and is discussed in more detail below, see \textbf{Section 4}. $%
\blacksquare $

The fact that the $8$ coefficients $c_{n\ell }$ ($n=1,2;$ $\ell =1,2,3,4$)
are expressed by the formulas (\ref{2}) in terms of $7$ \textit{a priori
arbitrary} parameters might suggest that these $8$ coefficients $c_{n\ell }$
are required to satisfy a \textit{single} condition. As it turns out, they
are in fact required to satisfy $2$ restrictions: and various \textit{%
explicit} avatars of these \textit{constraints}---\textit{sufficient} to
guarantee that the system (\ref{1}) possess the solution described by
\textbf{Proposition 1-1}---are provided below (see \textbf{Subsection 3.2}
and\textbf{\ Appendix C}). This is discussed below (after the proof of
\textbf{Proposition 1-1}:\textbf{\ }see \textbf{Section 2}) in the context
(see \textbf{Section 3}) of the related---quite important in \textit{%
applicative} contexts---issue of the \textit{inversion} of the relations (%
\ref{2}): namely the task of expressing in terms of the $8$ coefficients $%
c_{n\ell }$ ($n=1,2;$ $\ell =1,2,3,4$) the $7$ parameters $a_{n}$, $b_{n}$
and $\gamma _{j}$ ($n=1,2;$ $j=1,2,3$)---hence of obtaining the \textit{%
rather explicit} solution of the system (\ref{1}) provided by \textbf{%
Proposition 1.1}, as well as the\textit{\ }$2$ \textit{constraints} on the $%
8 $ coefficients $c_{n\ell }$ which entail this possibility.

Finally, let us suggest that readers primarily interested in the utilization
in \textit{applicative }contexts of the \textit{solvable} subclass of the
system of ODEs (\ref{1}) (and/or (\ref{1Iso})) have immediately a look below
at \textbf{Section\ 5}, to assess if, and how, the findings reported in this
paper are indeed likely to be useful for their purposes; and they might also
be interested to have a quick look at \textbf{Subsection 3.5}, focussed on
the subcase of the system (\ref{1}) with $c_{14}=c_{21}=0$, which is
particularly relevant in \textit{applicative} contexts.

\textbf{Remark 1-5}. Hereafter we assume that \textit{none} of the $4$
parameters $a_{1},$ $a_{2},$ $b_{1},$ $b_{2}$ vanishes; since it is plain
from \textbf{Proposition 1-1}---see in particular the eqs. (\ref{x1x2wt}%
)---that if anyone of these $4$ parameters vanishes the findings reported in
this paper become rather less interesting. $\blacksquare $

\section{Proof of Proposition 1-1}

The starting point of our treatment is the trivially solvable ODE
\begin{equation}
\dot{y}=y^{3}~,  \label{ydot}
\end{equation}%
the solution of which---with \textit{arbitrary} initial-value $y\left(
0\right) $---is easily seen to coincide with (\ref{Solyt}).

Next, we introduce the ODE
\begin{equation}
\dot{w}=w^{3}+\gamma _{1}yw^{2}+\gamma _{2}y^{2}w+\gamma _{3}y^{3}~,
\label{wdot}
\end{equation}%
and we set (see (\ref{wwwt}))
\begin{subequations}
\begin{equation}
w\left( t\right) =u\left( t\right) y\left( t\right) ~,  \label{wu}
\end{equation}%
implying of course
\begin{equation}
u\left( t\right) =w\left( t\right) /y\left( t\right) ~,
\end{equation}%
hence%
\begin{equation}
\dot{w}\left( t\right) =\dot{u}\left( t\right) y\left( t\right) +u\left(
t\right) \dot{y}\left( t\right) ~,
\end{equation}%
and (via (\ref{ydot}))%
\begin{eqnarray}
&&\dot{u}=y^{2}\left[ u^{3}+\gamma _{1}u^{2}+\left( \gamma _{2}-1\right)
u+\gamma _{3}\right] ~,  \nonumber \\
&&\dot{u}\left[ u^{3}+\gamma _{1}u^{2}+\left( \gamma _{2}-1\right) u+\gamma
_{3}\right] ^{-1}=y^{2}~,  \nonumber \\
&&\dot{u}\left( \frac{\lambda _{1}}{u-u_{1}}+\frac{\lambda _{2}}{u-u_{2}}+%
\frac{\lambda _{3}}{u-u_{3}}\right) =\left[ y\left( 0\right) \right]
^{2}\left\{ 1-2\left[ y\left( 0\right) \right] ^{2}t\right\} ^{-1}~,
\label{2ODEut}
\end{eqnarray}%
where the $3$ parameters $\gamma _{j}$ are \textit{a priori arbitrary }and
the $6$ parameters $\lambda _{j}$ and $u_{j}$ are clearly defined in terms
of them by the relations (\ref{alphauj}) and (\ref{landuj}).

The formula (\ref{1ut}) is then an immediate consequence of this ODE (\ref%
{2ODEut}).

Finally, it is easy to verify---via a straightforward if tedious
computation---that the relations (\ref{x12t}), implying of course
\end{subequations}
\begin{subequations}
\begin{equation}
\dot{x}_{1}\left( t\right) =\left[ b_{2}\dot{y}\left( t\right) -a_{2}\dot{w}%
\left( t\right) \right] /c~,~~~\dot{x}_{2}\left( t\right) =\left[ -b_{1}\dot{%
y}\left( t\right) +a_{1}\dot{w}\left( t\right) \right] /c~,
\end{equation}%
hence, via (\ref{ydot}) and (\ref{wdot}),%
\begin{eqnarray}
\dot{x}_{1} &=&\left[ b_{2}y^{3}-a_{2}\left( w^{3}+\gamma _{1}w^{2}y+\gamma
_{2}wy^{2}+\gamma _{3}y^{3}\right) \right] /c~,  \nonumber \\
\dot{x}_{2} &=&\left[ -b_{1}y^{3}+a_{1}\left( w^{3}+\gamma _{1}w^{2}y+\gamma
_{2}wy^{2}+\gamma _{3}y^{3}\right) \right] /c~,
\end{eqnarray}%
coincide---via the relations (\ref{ywxt})---with the system (\ref{1}),
provided the $8$ parameters $c_{n\ell }$ in the right-hand sides of (\ref{1}%
) are given by the formulas (\ref{2}) (with (\ref{1KK})) in terms of the $7$
parameters $\gamma _{1},$ $\gamma _{2},$ $\gamma _{3},$ $a_{1},$ $a_{2},$ $%
b_{1},$ $b_{2}$.

\textbf{Proposition 1-1} is thereby proven.

\section{Inversion of the relations (\protect\ref{2}) and constraints on the
parameters $c_{n\ell }$}

The task accomplished in this \textbf{Section 3} is nontrivial, due to the
\textit{nonlinear} character of the algebraic equations to be solved: it
could not be performed by the blind use of computer manipulations programs
such as \textbf{Mathematica}; and it is clear that there is not a \textit{%
unique} way to manage it.

The results reported in this \textbf{Section 3} are likely to be
particularly relevant for researchers who are interested to use in \textit{%
applicative} contexts the findings reported in this paper.

\textbf{Remark 3-1}. Note that below we often take advantage of the
definition (\ref{cc}) of $c$. $\blacksquare $

\subsection{First step towards inverting the relations (\protect\ref{2})}

The following easy consequences of the relations (\ref{2}) constitute a
first step towards their inversion.

From (\ref{1c11}) and (\ref{1c21}),
\end{subequations}
\begin{subequations}
\label{31aacc}
\begin{equation}
a_{1}c_{11}+a_{2}c_{21}=\left( a_{1}\right) ^{3}~;  \label{ac11c21}
\end{equation}%
from (\ref{1c12}) and (\ref{1c22}),%
\begin{equation}
a_{1}c_{12}+a_{2}c_{22}=3\left( a_{1}\right) ^{2}a_{2}~;  \label{ac12c22}
\end{equation}%
from (\ref{1c13}) and (\ref{1c23}),%
\begin{equation}
a_{1}c_{13}+a_{2}c_{23}=3a_{1}\left( a_{2}\right) ^{2}~;  \label{ac13c23}
\end{equation}%
from (\ref{1c14}) and (\ref{1c24}),%
\begin{equation}
a_{1}c_{14}+a_{2}c_{24}=\left( a_{2}\right) ^{3}~.  \label{ac14c24}
\end{equation}

\textbf{Remark 3.1-1}. Clearly these $4$ algebraic equations imply the
following $2$ neat consequences:
\end{subequations}
\begin{equation}
\sum_{n=1}^{2}\left[ a_{n}\left( c_{n1}+sc_{n2}+c_{n3}+sc_{n4}\right) \right]
=\left( a_{1}+sa_{2}\right) ^{3}~,~~~s=\pm ~.
\end{equation}%
But we shall not make use of these particular relation below. $\blacksquare $

\subsection{Second step: determination of the $2$ parameters $a_{n}$ in
terms of the $8$ coefficients $c_{n\ell }$ and of $2$ constraints on these $%
8 $ coefficients}

It is plain that the $4$ equations (\ref{31aacc}) featuring the $10$
dependent variables $a_{n}$ and $c_{n\ell }$ ($n=1,2$; $\ell =1,$ $2$, $3$, $%
4$) are independent of each other (since they involve different variables);
and since they involve only the $2$ parameters $a_{1}$ and $a_{2}$ in
addition to the $8$ coefficients $c_{n\ell }$ it is quite natural to infer
that they determine these $2$ parameters $a_{1}$ and $a_{2}$ in terms of the
$8$ coefficients $c_{n\ell }$ and that they moreover imply that the $8$
coefficients $c_{n\ell }$ satisfy $2$ \textit{constraints}. In this \textbf{%
Subsection 3.2} we indeed obtain various explicit expressions of the $2$
parameters $a_{n}$ in terms of the $8$\textit{\ coefficients }$c_{n\ell }$
and we also obtain various \textit{avatars} of $2$\textit{\ constraints}
implied by the $4$ eqs.(\ref{31aacc}) for these $8$ coefficients (with more
implied by the formulas reported in \textbf{Appendices A}, \textbf{B} and
\textbf{C}).

The determination from the $4$ eqs. (\ref{31aacc}) of the $2$ parameters $%
a_{1}$ and $a_{2}$ in terms of the $8$ coefficients $c_{n\ell }$, and of $2$
constraints on the $8$ coefficients $c_{n\ell }$, can be achieved via
different routes. A convenient way to make some initial progress is by
introducing the auxiliary quantity
\begin{equation}
\alpha =a_{2}/a_{1}~.  \label{32alfa}
\end{equation}%
Then, by taking the ratios of each of the (last $3$ of the $4$) equations (%
\ref{31aacc}) over that written before it, one easily gets the following $3$
\textit{quadratic} equations for the quantity $\alpha $:
\begin{subequations}
\label{32alfa2}
\begin{equation}
3c_{21}\alpha ^{2}+(3c_{11}-c_{22})\alpha -c_{12}=0~,  \label{32alfa2a}
\end{equation}

\begin{equation}
c_{22}\alpha ^{2}+(c_{12}-c_{23})\alpha -c_{13}=0~,  \label{32alfa2b}
\end{equation}

\begin{equation}
c_{23}\alpha ^{2}+(c_{13}-3c_{24})\alpha -3c_{14}=0~.  \label{32alfa2c}
\end{equation}

Next, multiply (\ref{32alfa2b}) by $3c_{21}$ and subtract from it (\ref%
{32alfa2a}) itself multiplied by $c_{22}$; likewise multiply (\ref{32alfa2c}%
) by $c_{22}$ and subtract from it (\ref{32alfa2b}) itself multiplied by $%
c_{23}$; and finally multiply (\ref{32alfa2a}) by $c_{23}$ and subtract from
it (\ref{32alfa2c}) itself multiplied by $3c_{21}$: in this manner the
following $3$ explicit expressions of $\alpha $ in terms of the $8$
coefficients $c_{n\ell }$ are easily obtained:
\end{subequations}
\begin{subequations}
\label{32alfa}
\begin{equation}
\alpha =\frac{c_{12}c_{22}-3c_{13}c_{21}}{%
3(c_{11}c_{22}-c_{12}c_{21}+c_{21}c_{23})-(c_{22})^{2}}~,  \label{32alfaa}
\end{equation}%
\begin{equation}
\alpha =\frac{c_{13}c_{23}-3c_{14}c_{22}}{%
c_{12}c_{23}-c_{13}c_{22}-(c_{23})^{2}+3c_{22}c_{24}}~,  \label{32alfab}
\end{equation}%
\begin{equation}
\alpha =\frac{c_{12}c_{23}-9c_{14}c_{21}}{%
3(c_{11}c_{23}-c_{13}c_{21})+9c_{21}c_{24}-c_{22}c_{23}}~.  \label{32alfac}
\end{equation}

\textbf{Remark 3.2-1}. Of course these are only $3$ specific expressions of $%
\alpha $, arbitrarily selected out of a plurality of \textit{different}---
but as well \textit{valid}---formulas which might be obtained by analogous
developments. And note that this remark is as well relevant for several of
the following developments. $\blacksquare $

The simultaneous validity of these $3$ expressions of $\alpha $ implies of
course the following $3$ relations among the parameters $c_{n\ell }$:
\end{subequations}
\begin{subequations}
\label{32Constraints}
\begin{eqnarray}
&&\left( c_{12}c_{22}-3c_{13}c_{21}\right) \left[ c_{12}c_{23}-c_{13}c_{22}-%
\left( c_{23}\right) ^{2}+3c_{22}c_{24}\right]  \nonumber \\
&=&\left( c_{13}c_{23}-3c_{14}c_{22}\right) \left[ 3\left(
c_{11}c_{22}-c_{12}c_{21}+c_{21}c_{23}\right) -\left( c_{22}\right) ^{2}%
\right] ~,  \label{32Constraintsa}
\end{eqnarray}%
\begin{eqnarray}
&&\left( c_{13}c_{23}-3c_{14}c_{22}\right) \left[ 3\left(
c_{11}c_{23}-c_{13}c_{21}\right) +9c_{21}c_{24}-c_{22}c_{23}\right]  \nonumber
\\
&=&\left( c_{12}c_{23}-9c_{14}c_{21}\right) \left[
c_{12}c_{23}-c_{13}c_{22}-\left( c_{23}\right) ^{2}+3c_{22}c_{24}\right] ~,
\label{32Constraintsb}
\end{eqnarray}%
\begin{eqnarray}
&&\left( c_{12}c_{23}-9c_{14}c_{21}\right) \left[ 3\left(
c_{11}c_{22}-c_{12}c_{21}+c_{21}c_{23}\right) -\left( c_{22}\right) ^{2}%
\right]  \nonumber \\
&=&\left( c_{12}c_{22}-3c_{13}c_{21}\right) \left[ 3\left(
c_{11}c_{23}-c_{13}c_{21}\right) +9c_{21}c_{24}-c_{22}c_{23}\right] ~.
\label{32Constraintsc}
\end{eqnarray}%
Clearly---as implied by their derivation from the relations (\ref{32alfa}%
)---any \textit{two} of these $3$ relations imply the \textit{third}; so
they provide at most $2$ \textit{independent constraints} on the $8$
coefficients $c_{n\ell }$. But in fact they only entail $1$ \textit{%
constraint} on the $8$ coefficients $c_{n\ell }$, as demonstrated by the
following remarkable phenomenon: if \textit{any one} of the $3$ eqs. (\ref%
{32Constraints}) is \textit{solved} for any specific one of the $8$
coefficients $c_{n\ell },$ then the $3$ results thereby obtained from these $%
3$ equations---in terms of the other $7$ coefficients $c_{n\ell }$---are
\textit{identical}. Since each of the $3$ eqs. (\ref{32Constraints}) is
\textit{linear} or \textit{quadratic} in each of the $8$ coefficients $%
c_{n\ell }$, these operations can be explicitly performed by hand (although
it is wise to also check the result via \textbf{Mathematica}). The
corresponding formulas are reported in \textbf{Appendix A}. Each of the $8$
formulas reported there provides the \textit{explicit} expression of one of
the $8$ parameters $c_{n\ell }$ in term of the other $7$; so each of them
provides a \textit{constraint} on the $8$ coefficients $c_{n\ell }$; and a
\textit{constraint} is as well provided by each of the $3$ eqs. (\ref%
{32Constraints}), or by any reasonable combinations of all these formulas.
Of course \textit{all} these constraints are \textit{equivalent}.
Nevertheless their \textit{explicit} exhibition---especially in the "solved"
version detailed in \textbf{Appendix A}---is worthwhile because of its
potential usefulness in applicative contexts.

Clearly any one of the $3$ formulas (\ref{32alfa}) provides an expression of
$\alpha $ in terms of the coefficients $c_{n\ell }$; these $3$ expressions
of $\alpha $ are of course \textit{equivalent} if the coefficients $c_{n\ell
}$ satisfy the \textit{constraint} mentioned above (see (\ref{32Constraints}%
) and \textbf{Appendix A}), as we hereafter assume. It is then an easy task
to get various \textit{explicit} expressions of the parameters $a_{1}$ and $%
a_{2}$ from the equations (\ref{31aacc}); for instance, by dividing (\ref%
{ac11c21}) by $a_{1}$ one gets
\end{subequations}
\begin{subequations}
\label{32aa12}
\begin{equation}
\left( a_{1}\right) ^{2}=c_{11}+\alpha c_{21}~,  \label{32a1}
\end{equation}%
and then, from (\ref{ac12c22}),%
\begin{equation}
a_{2}=\frac{a_{1}c_{12}}{3\left( a_{1}\right) ^{2}-c_{22}}=\frac{a_{1}c_{12}%
}{3\left( c_{11}+\alpha c_{21}\right) -c_{22}}~,  \label{32a2}
\end{equation}%
where the second equality is of course implied by (\ref{32a1}) (and of
course in these formulas $\alpha $ is given by any one of the $3$ formulas (%
\ref{32alfa})).

There still remains the task to determine the \textit{second constraint} on
the $8$ coefficients implied by the $4$ eqs. (\ref{31aacc}).\ Other explicit
expressions of the $2$ parameters $a_{1}$ and $a_{2}$---or rather of their
squares---in terms of the $8$ coefficients $c_{n\ell }$ are also provided
below in this context.

Let us begin by performing the following $4$ operations on the $4$ eqs. (\ref%
{31aacc}). ($1$) We multiply the \textit{first} of these $4$ equations by $%
c_{22}$ and subtract from the result the \textit{second} of these $4$
equations itself multiplied by $c_{21}$. ($2$) We multiply the \textit{third}
of these $4$ equations by $c_{14}$ and subtract from the result the \textit{%
fourth} of these $4$ equations itself multiplied by $c_{13}$. ($3$) We
multiply the \textit{first} of these $4$ equations by $c_{23}$ and subtract
from the result the \textit{third} of these $4$ equations itself multiplied
by $c_{21}$. ($4$) We multiply the \textit{second} of these $4$ equations by
$c_{14}$ and subtract from the result the \textit{fourth} of these $4$
equations itself multiplied by $c_{12}$. There thus obtain the following $4$
equations:
\end{subequations}
\begin{subequations}
\label{324eqs}
\begin{equation}
c_{11}c_{22}-c_{12}c_{21}=c_{22}\left( a_{1}\right) ^{2}-3c_{21}a_{1}a_{2}~,
\end{equation}%
\begin{equation}
c_{14}c_{23}-c_{13}c_{24}=-c_{13}\left( a_{2}\right) ^{2}+3c_{14}a_{1}a_{2}~,
\end{equation}%
\begin{equation}
c_{11}c_{23}-c_{13}c_{21}=c_{23}\left( a_{1}\right) ^{2}-3c_{21}\left(
a_{2}\right) ^{2}~,
\end{equation}%
\begin{equation}
c_{14}c_{22}-c_{12}c_{24}=-c_{12}\left( a_{2}\right) ^{2}+3c_{14}\left(
a_{1}\right) ^{2}~.
\end{equation}

Next, we perform the following cycle of analogous operations on these eqs. (%
\ref{324eqs}) (which are clearly equivalent to the $4$ eqs. (\ref{31aacc})).
($1$) We multiply the \textit{first} of these $4$ equations by $c_{14}$ and
we add the result to the \textit{second} of these $4$ equations itself
multiplied by $c_{21}.$ ($2$) We multiply the \textit{third} of these $4$
equations by $c_{12}$ and we subtract from the result the \textit{fourth} of
these $4$ equations itself multiplied by $3c_{21}$. ($3$) We multiply the
\textit{third} of these $4$ equations by $3c_{14}$ and we subtract from the
result the \textit{fourth} of these $4$ equations itself multiplied by $%
c_{23}$. There thus obtain the following $3$ equations:
\end{subequations}
\begin{subequations}
\label{32cc1234}
\begin{equation}
c_{14}\left( c_{11}c_{22}-c_{12}c_{21}\right) +c_{21}\left(
c_{14}c_{23}-c_{13}c_{24}\right) =c_{14}c_{22}\left( a_{1}\right)
^{2}-c_{13}c_{21}\left( a_{2}\right) ^{2}~,  \label{32cc1234a}
\end{equation}%
\begin{equation}
c_{12}\left( c_{11}c_{23}-c_{13}c_{21}\right) -3c_{21}\left(
c_{14}c_{22}-c_{12}c_{24}\right) =\left( c_{12}c_{23}-9c_{14}c_{21}\right)
\left( a_{1}\right) ^{2}~,  \label{32cc1234b}
\end{equation}%
\begin{equation}
3c_{14}\left( c_{11}c_{23}-c_{13}c_{21}\right) -c_{23}\left(
c_{14}c_{22}-c_{12}c_{24}\right) =\left( c_{12}c_{23}-9c_{14}c_{21}\right)
\left( a_{2}\right) ^{2}~.  \label{32cc1234c}
\end{equation}

Next, we multiply the second of these $3$ equations (\ref{32cc1234}) by $%
c_{14}c_{22}$ and we subtract from the result the last of these $3$
equations itself multiplied by $c_{13}c_{21}$.We thereby get the following
equation:
\end{subequations}
\begin{eqnarray}
&&c_{14}\left( c_{11}c_{23}-c_{13}c_{21}\right) \left(
c_{12}c_{22}-3c_{13}c_{21}\right)  \nonumber \\
&&+c_{21}\left( c_{14}c_{22}-c_{12}c_{24}\right) \left(
c_{13}c_{23}-3c_{14}c_{22}\right)  \nonumber \\
&=&\left( c_{12}c_{23}-9c_{14}c_{21}\right) \left[ c_{14}c_{22}\left(
a_{1}\right) ^{2}-c_{13}c_{21}\left( a_{2}\right) ^{2}\right] ~,
\end{eqnarray}%
hence, by comparing the right-hand side of this equation with the right-hand
side of eq. (\ref{32cc1234a}), we get finally the following \textit{%
constraint} involving the $8$ coefficients $c_{n\ell }$:
\begin{eqnarray}
&&c_{14}\left( c_{11}c_{23}-c_{13}c_{21}\right) \left(
c_{12}c_{22}-3c_{13}c_{21}\right)  \nonumber \\
&&+c_{21}\left( c_{14}c_{22}-c_{12}c_{24}\right) \left(
c_{13}c_{23}-3c_{14}c_{22}\right)  \nonumber \\
&=&\left( c_{12}c_{23}-9c_{14}c_{21}\right) \left[ c_{14}\left(
c_{11}c_{22}-c_{12}c_{21}\right) +c_{21}\left(
c_{14}c_{23}-c_{13}c_{24}\right) \right] ~.  \label{32Constraint2}
\end{eqnarray}

\textbf{Remark 3.2-1}. This \textit{constraint} is invariant under the
transformation (\ref{1Change}). $\blacksquare $

But, as it happens, this is rather a \textit{new avatar} of the \textit{%
constraint} on the $8$ coefficients $c_{n\ell }$ already obtained above: it
is again a \textit{homogeneous polynomial of fifth degree} involving these $%
8 $ coefficients, but it is only \textit{linear} in the $2$ coefficients $%
c_{11}$ and $c_{24},$ and \textit{quadratic} in each of the other $6$
coefficients. Indeed, it is easily seen that this \textit{constraint}
implies the following expression of the coefficient $c_{11}$ in terms of the
other $7$ coefficients (and of course an analogous expression of the
coefficient $c_{24}$ in terms of the other $7$ coefficients can then be
obtained via the transformation (\ref{1Change})):
\begin{eqnarray}
c_{11} &=&\left\{ -3c_{12}c_{22}c_{24}+3c_{14}\left[ 3c_{12}c_{21}+\left(
c_{22}\right) ^{2}\right] \right.  \nonumber \\
&&-\left[ \left( c_{12}\right) ^{2}+9c_{14}c_{21}\right] c_{23}+c_{12}\left(
c_{23}\right) ^{2}  \nonumber \\
&&+c_{13}\left[ c_{12}c_{22}-c_{22}c_{23}+9c_{21}c_{24}\right]  \nonumber \\
&&\left. -3\left( c_{13}\right) ^{2}c_{21}\right\} /\left[
3(3c_{14}c_{22}-c_{13}c_{23})\right] ~.
\end{eqnarray}%
And it is easily seen that this formula is in fact \textit{identical} to eq.
(\ref{A1}).

And the same happens for the expressions of the other $7$ coefficients $%
c_{n\ell }$ obtained in an analogous manner: see \textbf{Appendix A}.

Finally, to make progress towards the identification of the \textit{second
constraint} we note that the $2$ eqs. (\ref{32cc1234b}) and (\ref{32cc1234c}%
) provide the following $2$ expressions of (the squares of) $a_{1}$ and $%
a_{2}$ in terms of the $8$ coefficients $c_{n\ell }$:
\begin{subequations}
\label{BAn}
\begin{equation}
\left( a_{n}\right) ^{2}=A_{n}\left( C\right) ~,~~~n=1,2~,  \label{32ansq}
\end{equation}%
\begin{equation}
A_{1}\left( C\right) \equiv \left[ c_{12}\left(
c_{11}c_{23}-c_{13}c_{21}\right) -3c_{21}\left(
c_{14}c_{22}-c_{12}c_{24}\right) \right] /\left(
c_{12}c_{23}-9c_{14}c_{21}\right) ~,
\end{equation}%
\begin{equation}
A_{2}\left( C\right) \equiv \left[ 3c_{14}\left(
c_{11}c_{23}-c_{13}c_{21}\right) -c_{23}\left(
c_{14}c_{22}-c_{12}c_{24}\right) \right] /\left(
c_{12}c_{23}-9c_{14}c_{21}\right) ~.
\end{equation}%
Here and below we denote by $C$ the set of $8$ coefficients $c_{n\ell }$.

Let us then return to the original $4$ eqs. (\ref{31aacc}), and let us
replace---in an obvious manner---the \textit{squares} of the variables $%
a_{n} $ by their expressions (\ref{32ansq}). We thus obtain the following
system of $4$ \textit{homogeneous linear} equations for the $2$ dependent
variables $a_{n}$:
\end{subequations}
\begin{subequations}
\begin{equation}
\left[ c_{11}-A_{1}\left( C\right) \right] a_{1}+c_{21}a_{2}=0~,
\end{equation}%
\begin{equation}
c_{12}a_{1}+\left[ c_{22}-3A_{1}\left( C\right) \right] a_{2}=0~,
\end{equation}%
\begin{equation}
\left[ c_{13}-3A_{2}\left( C\right) \right] a_{1}+c_{23}a_{2}=0~,
\end{equation}%
\begin{equation}
c_{14}a_{1}+\left[ c_{24}-A_{2}\left( C\right) \right] a_{2}=0~.
\end{equation}%
By selecting any pair of these $4$ equations we get $6$ different systems of
$2$ \textit{homogeneous linear} equations, which must be satisfied by the $2$
(nonvanishing!) dependent variables $a_{1}$ and $a_{2}.$ Hence the $6$
determinants of the coefficients of these $6$ systems must vanish, providing
thereby---in principle---as many \textit{constraints} on the $8$
coefficients $c_{n\ell }$. But $2$ of the resulting relations yield an
\textit{identically vanishing} result ($0=0$), and the other $4$ the \textit{%
same} outcome (provided the $8$ coefficients have \textit{generic} values,
none of them vanishing): yielding the following \textit{second constraint}
on the $8$ coefficients $c_{n\ell }$:
\end{subequations}
\begin{eqnarray}
&&\left[ c_{12}\left( c_{13}-3c_{24}\right) +3c_{14}\left(
-3c_{11}+c_{22}\right) \right] \left[ 3c_{21}\left( c_{13}-3c_{24}\right)
+c_{23}\left( c_{22}-3c_{11}\right) \right]  \nonumber \\
&=&(c_{12}c_{23}-9c_{14}c_{21})^{2}~.  \label{32SecCona}
\end{eqnarray}%
The fact that this \textit{second constraint} is \textit{not equivalent} to
the \textit{first} is demonstrated by solving it for each of the $8$
coefficients $c_{n\ell }$ and by noting that the resulting expressions of
the $8$ coefficients $c_{n\ell }$---as reported in \textbf{Appendix B}---do
\textit{not} coincide with those reported in \textbf{Appendix A}.

The formulas reported in \textbf{Appendices A} and \textbf{B} become of
course pairwise \textit{equivalent} if the $8$ coefficients $c_{n\ell }$ are
required to satisfy the $2$ \textit{constraints} identified in this \textbf{%
Subsection 3.2}: both of them!

The way is now open---solving simultaneously \textit{both} the $2$ \textit{%
independent constraints} satisfied by the $8$ coefficients $c_{n\ell }$---to
obtain explicit expressions of \textit{each pair} of the $8$ coefficients $%
c_{n\ell }$ in terms of the other $6$: remarkably, in all cases the relevant
formulas could be \textit{explicitly} obtained (via \textbf{Mathematica})
although in some cases they are too complicated to be usefully displayed:
see \textbf{Appendix C}. To overcome as much as possible this difficulty the
following remark turned out to be useful.

\textbf{Remark 3.2-2}. The \textit{second constraint} (\ref{32SecCona}) is a
\textit{homogeneous polynomial} equation of degree $4$ in the $8$ variables $%
c_{n\ell }$; and several of the \textit{equivalent} versions of the \textit{%
first constraint} are \textit{also homogeneous polynomial} equations, in
particular the $3$ versions (\ref{32Constraints}) are each \textit{%
homogeneous polynomial} equations of degree $4$; but $2$ of the equations
displayed in \textbf{Appendix A}---those containing no square-roots: see (%
\ref{A1}) and (\ref{A4})---are (\textit{de facto}) \textit{homogeneous
polynomial} equations of degree (only!) $3$. $\blacksquare $

\subsection{Third step: determination of the $2$ parameters $b_{n}$}

The starting point for the next steps to invert the $8$ relations (\ref{2})
are the following $4$ relations, obtained by summing pairwise the $8$
relations (\ref{2}) appropriately multiplied by $b_{1}$ and $b_{2}$:
\begin{subequations}
\begin{equation}
b_{1}c_{11}+b_{2}c_{21}=K_{1}~,
\end{equation}%
\begin{equation}
b_{1}c_{12}+b_{2}c_{22}=K_{2}~,
\end{equation}%
\begin{equation}
b_{1}c_{13}+b_{2}c_{23}=K_{3}~,
\end{equation}%
\begin{equation}
b_{1}c_{14}+b_{2}c_{24}=K_{4}~,
\end{equation}%
of course with the $4$ quantities $K_{\ell }$ defined by the identities (\ref%
{1KK}).

We now note that these $4$ relations can be re-written, via (\ref{1KK}), as
follows
\end{subequations}
\begin{subequations}
\label{33KK}
\begin{equation}
\gamma _{1}a_{1}\left( b_{1}\right) ^{2}+\gamma _{2}\left( a_{1}\right)
^{2}b_{1}+\gamma _{3}\left( a_{1}\right) ^{3}=-\left( b_{1}\right)
^{3}+b_{1}c_{11}+b_{2}c_{21}~,  \label{33K1}
\end{equation}%
\begin{eqnarray}
&&\gamma _{1}\left[ a_{2}\left( b_{1}\right) ^{2}+2a_{1}b_{1}b_{2}\right]
+\gamma _{2}\left[ \left( a_{1}\right) ^{2}b_{2}+2a_{1}a_{2}b_{1}\right]
+3\gamma _{3}\left( a_{1}\right) ^{2}a_{2}  \nonumber \\
&=&-3\left( b_{1}\right) ^{2}b_{2}+b_{1}c_{12}+b_{2}c_{22}~,  \label{33K2}
\end{eqnarray}%
\begin{eqnarray}
&&\gamma _{1}\left[ a_{1}\left( b_{2}\right) ^{2}+2a_{2}b_{1}b_{2}\right]
+\gamma _{2}\left[ \left( a_{2}\right) ^{2}b_{1}+2a_{1}a_{2}b_{2}\right]
+3\gamma _{3}a_{1}\left( a_{2}\right) ^{2}  \nonumber \\
&=&-3b_{1}\left( b_{2}\right) ^{2}+b_{1}c_{13}+b_{2}c_{23}~,  \label{33K3}
\end{eqnarray}%
\begin{equation}
\gamma _{1}a_{2}\left( b_{2}\right) ^{2}+\gamma _{2}\left( a_{2}\right)
^{2}b_{2}+\gamma _{3}\left( a_{2}\right) ^{3}=-\left( b_{2}\right)
^{3}+b_{1}c_{14}+b_{2}c_{24}~.  \label{33K4}
\end{equation}

These are clearly $4$ \textit{linear} equations satisfied by the $3$
parameters $\gamma _{j}$ ($j=1,2,3$); so that these $3$ parameters $\gamma
_{j}$ can then be \textit{explicitly} determined by solving the system
formed by any $3$ of these $4$ equations, obtaining thereby expressions of
these $3$ parameters $\gamma _{j}$ in terms of the $8$ coefficients $%
c_{n\ell }$ appearing in the right-hand sides of these equations (\ref{33KK}%
), of the $2$ parameters $a_{n}$---themselves already determined in terms of
the $8$ coefficients $c_{n\ell }$ (see \textbf{Subsection 3.2})---and of the
still undetermined $2$ parameters $b_{n}$. Once this step has been
completed, the unused one of the $4$ equations (\ref{33KK}) provides a
\textit{single explicit constraint} on the $2$ parameters $b_{n},$ implying
the determination of one of them in terms of the other, or of their ratio
(see below). But this rather natural route leads to quite complicated final
equations$.$ An equivalent, more practical, route is described below.

\textbf{Remark 3.3-1}. It is plain---by summing these $4$ relations (\ref%
{33KK}), with the second and fourth multiplied by $s=\pm $---that they imply
the $2$ relations
\end{subequations}
\begin{eqnarray}
&&\gamma _{1}\left( a_{1}+sa_{2}\right) \left( b_{1}+sb_{2}\right)
^{2}+\gamma _{2}\left( a_{1}+sa_{2}\right) ^{2}\left( b_{1}+sb_{2}\right)
+\gamma _{3}\left( a_{1}+sa_{2}\right) ^{3}  \nonumber \\
&=&-\left( b_{1}+sb_{2}\right) ^{3}+b_{1}\left(
c_{11}+sc_{12}+c_{13}+sc_{14}\right) +b_{2}\left(
c_{21}+sc_{22}+c_{23}+sc_{24}\right) ~,  \nonumber \\
&&s=\pm ~.
\end{eqnarray}%
But again we shall not need to take advantage of this relation. $%
\blacksquare $

\textbf{Remark 3.3-2}. By multiplying the formula (\ref{33K1}) by $%
(a_{2})^{3}$ and by then subtracting from the result the formula (\ref{33K4}%
) itself multiplied by $(a_{1})^{3}$ one gets (after some easy
simplifications) the formula
\begin{subequations}
\label{33gammas}
\begin{eqnarray}
&&a_{1}a_{2}\left( a_{1}b_{2}-a_{2}b_{1}\right) \left[ \left(
a_{1}b_{2}+a_{2}b_{1}\right) \gamma _{1}+a_{1}a_{2}\gamma _{2}\right]  \nonumber
\\
&=&-\left( a_{1}b_{2}\right) ^{3}+\left( a_{2}b_{1}\right) ^{3}  \nonumber \\
&&+\left( a_{1}\right) ^{3}\left( b_{1}c_{14}+b_{2}c_{24}\right) -\left(
a_{2}\right) ^{3}\left( b_{1}c_{11}+b_{2}c_{21}\right) ;
\end{eqnarray}%
likewise, by multiplying the eq. (\ref{33K2}) by $a_{2}$ and by then
subtracting from the result the formula (\ref{33K3}) itself multiplied by $%
a_{1}$ one gets (after some easy simplifications) the formula%
\begin{eqnarray}
&&\left( a_{1}b_{2}-a_{2}b_{1}\right) \left[ \left(
a_{1}b_{2}+a_{2}b_{1}\right) \gamma _{1}+a_{1}a_{2}\gamma _{2}\right]  \nonumber
\\
&=&-3\left( a_{1}b_{2}-a_{2}b_{1}\right) b_{1}b_{2}+b_{2}\left(
a_{1}c_{23}-a_{2}c_{22}\right)  \nonumber \\
&&-b_{1}\left( a_{1}c_{13}-a_{2}c_{12}\right) ~.
\end{eqnarray}

Next, we multiply the second of these $2$ formulas by $a_{1}a_{2}$ and then
subtract the result from the first of these $2$ formulas, getting thereby
(after some easy simplifications) the formula%
\begin{eqnarray}
&&\left( a_{1}b_{2}-a_{2}b_{1}\right) ^{3}=b_{2}\left[ \left( a_{1}\right)
^{3}c_{24}-\left( a_{2}\right) ^{3}c_{21}-a_{1}a_{2}\left(
a_{1}c_{23}-a_{2}c_{22}\right) \right]  \nonumber \\
&&+b_{1}\left[ \left( a_{1}\right) ^{3}c_{14}-\left( a_{2}\right)
^{3}c_{11}-a_{1}a_{2}\left( a_{1}c_{13}-a_{2}c_{12}\right) \right] ~.
\label{33bb}
\end{eqnarray}

Finally it is convenient to re-write this equation via the following
positions:
\end{subequations}
\begin{equation}
a_{2}=\alpha a_{1}~,~~~c_{n\ell }=\eta \hat{c}_{n\ell }~,~~~b_{1}=\eta
~,~~~b_{2}=\alpha \beta \eta ~.  \label{betaeta}
\end{equation}%
Here of course $\alpha $ is the parameter already introduced in \textbf{%
Subsection 3.2} (see (\ref{32alfa})) and determined there in terms of the $8$
coefficients $c_{n\ell }$; $\eta $ is a parameter characterizing a rescaling
of the $8$ coefficients $c_{n\ell }$---introduced here to make notational
contact with the invariance property of the system (\ref{1}) mentioned in
\textbf{Remark 1-1} (see (\ref{1Inveta})) which implies that this parameter
can be eventually altogether eliminated---i. e., assigned an \textit{%
arbitrary} \textit{nonvanishing value }(for instance, just the value $\eta
=1 $)---via a corresponding appropriate rescaling of the independent
variable $t $; $\beta $ is the parameter that we determine immediately below
in terms of the $8$ coefficients $c_{n\ell }$; and the last $2$ eqs. (\ref%
{betaeta}) determine of course the $2$ parameters $b_{1}$ and $b_{2}$,
thereby completing the task indicated by the title of this \textbf{%
Subsection 3.3}.

To determine the parameter $\beta $ in terms of the $8$ parameters $c_{n\ell
}$---or, equivalently, $\hat{c}_{n\ell }$ (see (\ref{betaeta}) and (\ref%
{1Inveta}))---we insert in eq. (\ref{33bb}) the positions (\ref{betaeta}),
getting thereby the following \textit{cubic} equation for this parameter:
\begin{eqnarray}
&&\left( \beta -1\right) ^{3}=\beta \left( \hat{c}_{24}\alpha ^{-2}-\hat{c}%
_{23}\alpha ^{-1}+\hat{c}_{22}-\hat{c}_{21}\alpha ^{-1}\right)  \nonumber \\
&&+\hat{c}_{14}\alpha ^{-2}-\hat{c}_{13}\alpha ^{-1}+\hat{c}_{12}-\hat{c}%
_{11}\alpha ^{-1}~.  \label{33beta}
\end{eqnarray}%
\textit{Explicit} solutions of this equation can of course be provided via
the standard Cardano formulas.

\subsection{Fourth step: determination of the $3$ parameters $\protect\gamma %
_{j}$ ($j=1,2,3$) in terms of the $8$ coefficients $c_{n\ell }$}

The \textit{explicit} determination of the $3$ parameters $\gamma _{j}$ ($%
j=1,2,3$) is now in principle a standard task, amounting---as shown
below---to the solution of any triad of the several \textit{linear}
equations satisfied by these quantities: see for instance the $4$ eqs. (\ref%
{33KK}) and the $3$ eqs. (\ref{33gammas}), or appropriate linear
combinations of these equations. Of course care must be taken to use $3$
equations which are \textit{independent} of each other. The expressions
obtained in this manner in terms of the $8$ coefficients $c_{n\ell }$---and
of the parameters $a_{n}$ and $b_{n}$ themselves expressed in terms of the $%
8 $ coefficients $c_{n\ell }$ as explained in \textbf{Subsections 3.2} and
\textbf{3.3}---are of course only valid provided the $8$ coefficients $%
c_{n\ell }$ satisfy the explicit \textit{constraints} determined above (see
\textbf{Subsection 3.2} and \textbf{Appendices A}, \textbf{B }and \textbf{C}%
); and they are \textit{equivalent}, but they need not look \textit{identical%
}.

The most straightforward procedure is to use $3$ out of the $4$ eqs. (\ref%
{33KK}): there are then of course $4$ possible choices.

\textbf{Remark 3.4-1}. Hereafter $c=a_{1}b_{2}-a_{2}b_{1},$ see (\ref{cc}). $%
\blacksquare $

The first choice are the $3$ eqs. (\ref{33K1}), (\ref{33K2}) and (\ref{33K3}%
). One then gets (if need be, with the help of \textbf{Mathematica)}
\begin{subequations}
\label{331gamj}
\begin{eqnarray}
&&\gamma _{1}=\left[ -3\left( a_{2}\right) ^{2}\left( b_{1}\right)
^{3}+6a_{1}a_{2}\left( b_{1}\right) ^{2}b_{2}-3\left( a_{1}\right)
^{2}b_{1}\left( b_{2}\right) ^{2}\right.  \nonumber \\
&&+3\left( a_{2}\right) ^{2}\left( b_{1}c_{11}+b_{2}c_{21}\right)
-2a_{1}a_{2}\left( b_{1}c_{12}+b_{2}c_{22}\right)  \nonumber \\
&&\left. +\left( a_{1}\right) ^{2}\left( b_{1}c_{13}+b_{2}c_{23}\right)
\right] /\left( a_{1}c^{2}\right) ~,  \label{331gam1}
\end{eqnarray}%
\begin{eqnarray}
&&\gamma _{2}=\left\{ 3\left( a_{2}\right) ^{2}\left( b_{1}\right)
^{4}-6a_{1}a_{2}\left( b_{1}\right) ^{3}b_{2}+3\left( a_{1}b_{1}b_{2}\right)
^{2}-3\left( a_{2}\right) ^{2}b_{1}\left( b_{1}c_{11}+b_{2}c_{21}\right)
\right.  \nonumber \\
&&\left. +3a_{1}a_{2}\left[ \left( b_{1}\right) ^{2}c_{12}-\left(
b_{2}\right) ^{2}c_{21}+b_{1}b_{2}\left( -c_{11}+c_{22}\right) \right]
\right.  \nonumber \\
&&\left. +\left( a_{1}\right) ^{2}\left[ -2\left( b_{1}\right)
^{2}c_{13}+\left( b_{2}\right) ^{2}c_{22}+b_{1}b_{2}\left(
c_{12}-2c_{23}\right) \right] \right\} /\left( a_{1}c\right) ^{2}~,
\label{331gam2}
\end{eqnarray}%
\begin{eqnarray}
&&\gamma _{3}=\left\{ -\left( a_{2}\right) ^{2}\left( b_{1}\right)
^{5}+2a_{1}a_{2}\left( b_{1}\right) ^{4}b_{2}-\left( a_{1}\right) ^{2}\left(
b_{1}\right) ^{3}\left( b_{2}\right) ^{2}\right.  \nonumber \\
&&+\left( a_{2}b_{1}\right) ^{2}\left( b_{1}c_{11}+b_{2}c_{21}\right)  \nonumber
\\
&&\left. -a_{1}a_{2}b_{1}\left[ \left( b_{1}\right) ^{2}c_{12}-\left(
b_{2}\right) ^{2}c_{21}+b_{1}b_{2}\left( -c_{11}+c_{22}\right) \right]
\right.  \nonumber \\
&&+\left( a_{1}\right) ^{2}\left[ \left( b_{1}\right) ^{3}c_{13}+\left(
b_{2}\right) ^{3}c_{21}+b_{1}\left( b_{2}\right) ^{2}\left(
c_{11}-c_{22}\right) \right.  \nonumber \\
&&\left. \left. -\left( b_{1}\right) ^{2}b_{2}\left( c_{12}-c_{23}\right)
\right] \right\} /\left[ \left( a_{1}\right) ^{3}c^{2}\right] ~.
\label{331gam3}
\end{eqnarray}

The second choice are the $3$ eqs. (\ref{33K1}), (\ref{33K2}) and (\ref{33K4}%
). One then gets (if need be, with the help of \textbf{Mathematica)}
\end{subequations}
\begin{subequations}
\label{332gamj}
\begin{eqnarray}
&&\gamma _{1}=\left[ -2\left( a_{2}b_{1}\right) ^{3}+3a_{1}\left(
a_{2}b_{1}\right) ^{2}b_{2}-\left( a_{1}b_{2}\right) ^{3}\right.  \nonumber \\
&&\left. +2\left( a_{2}\right) ^{3}\left( b_{1}c_{11}+b_{2}c_{21}\right)
-a_{1}\left( a_{2}\right) ^{2}\left( b_{1}c_{12}+b_{2}c_{22}\right) \right.
\nonumber \\
&&\left. +\left( a_{1}\right) ^{3}\left( b_{1}c_{14}+b_{2}c_{24}\right)
\right] /\left( a_{1}a_{2}c^{2}\right) ~,  \label{332gam1}
\end{eqnarray}%
\begin{eqnarray}
&&\gamma _{2}=\left\{ \left( a_{2}\right) ^{3}\left( b_{1}\right)
^{4}+2\left( a_{1}\right) ^{3}b_{1}\left( b_{2}\right) ^{3}-3\left(
a_{1}\right) ^{2}a_{2}\left( b_{1}b_{2}\right) ^{2}\right.  \nonumber \\
&&\left. -\left( a_{2}\right) ^{3}b_{1}\left( b_{1}c_{11}+b_{2}c_{21}\right)
+\left( a_{1}\right) ^{2}a_{2}b_{2}\left( b_{1}c_{12}+b_{2}c_{22}\right)
\right.  \nonumber \\
&&\left. +a_{1}\left( a_{2}\right) ^{2}\left[ \left( b_{1}\right)
^{2}c_{12}-3\left( b_{2}\right) ^{2}c_{21}+b_{1}b_{2}\left(
-3c_{11}+c_{22}\right) \right] \right.  \nonumber \\
&&\left. -2\left( a_{1}\right) ^{3}b_{1}\left(
b_{1}c_{14}+b_{2}c_{24}\right) \right\} /\left[ \left( a_{1}\right)
^{2}a_{2}c^{2}\right] ~,  \label{332gam2}
\end{eqnarray}%
\begin{eqnarray}
&&\gamma _{3}=\left\{ -\left( a_{2}\right) ^{2}\left( b_{1}\right)
^{4}b_{2}+2a_{1}a_{2}\left( b_{1}\right) ^{3}\left( b_{2}\right) ^{2}-\left(
a_{1}b_{1}\right) ^{2}\left( b_{2}\right) ^{3}\right.  \notag \\
&&+\left( a_{2}\right) ^{2}b_{1}b_{2}\left( b_{1}c_{11}+b_{2}c_{21}\right)
\nonumber \\
&&\left. +a_{1}a_{2}b_{2}\left[ -\left( b_{1}\right) ^{2}c_{12}+\left(
b_{2}\right) ^{2}c_{21}+b_{1}b_{2}\left( c_{11}-c_{22}\right) \right] \right.
\nonumber \\
&&\left. +\left( a_{1}b_{1}\right) ^{2}\left( b_{1}c_{14}+b_{2}c_{24}\right)
\right\} /\left[ \left( a_{1}\right) ^{2}a_{2}c^{2}\right] ~.
\label{332gam3}
\end{eqnarray}

The third choice are the $3$ eqs. (\ref{33K1}), (\ref{33K3}) and (\ref{33K4}%
). One then gets (if need be, with the help of \textbf{Mathematica)}
\end{subequations}
\begin{subequations}
\label{333gamj}
\begin{eqnarray}
&&\gamma _{1}=\left[ -\left( a_{2}b_{1}\right) ^{3}+3\left( a_{1}\right)
^{2}a_{2}b_{1}\left( b_{2}\right) ^{2}-2\left( a_{1}b_{2}\right) ^{3}\right.
\nonumber \\
&&+\left( a_{2}\right) ^{3}\left( b_{1}c_{11}+b_{2}c_{21}\right) -\left(
a_{1}\right) ^{2}a_{2}\left( b_{1}c_{13}+b_{2}c_{23}\right)  \nonumber \\
&&\left. +2\left( a_{1}\right) ^{3}\left( b_{1}c_{14}+b_{2}c_{24}\right)
\right] /\left( a_{1}a_{2}c^{2}\right) ~,  \label{333gam1}
\end{eqnarray}%
\begin{eqnarray}
&&\gamma _{2}=\left\{ 2\left( a_{2}b_{1}\right) ^{3}b_{2}-3a_{1}\left(
a_{2}b_{1}b_{2}\right) ^{2}+\left( a_{1}\right) ^{3}\left( b_{2}\right)
^{4}\right.  \nonumber \\
&&\left. +a_{1}\left( a_{2}\right) ^{2}b_{1}\left(
b_{1}c_{13}+b_{2}c_{23}\right) -2\left( a_{2}\right) ^{3}b_{2}\left(
b_{1}c_{11}+b_{2}c_{21}\right) \right.  \nonumber \\
&&\left. +\left( a_{1}\right) ^{2}a_{2}\left[ -3\left( b_{1}\right)
^{2}c_{14}+\left( b_{2}\right) ^{2}c_{23}+b_{1}b_{2}\left(
c_{13}-3c_{24}\right) \right] \right.  \nonumber \\
&&\left. -\left( a_{1}\right) ^{3}b_{2}\left( b_{1}c_{14}+b_{2}c_{24}\right)
\right\} /\left[ a_{1}\left( a_{2}\right) ^{2}c^{2}\right] ~,
\label{333gam2}
\end{eqnarray}%
\begin{eqnarray}
&&\gamma _{3}=\left\{ -\left( a_{2}\right) ^{2}\left( b_{1}\right)
^{3}\left( b_{2}\right) ^{2}-\left( a_{1}\right) ^{2}b_{1}\left(
b_{2}\right) ^{4}+2a_{1}a_{2}\left( b_{1}\right) ^{2}\left( b_{2}\right)
^{3}\right.  \nonumber \\
&&+\left( a_{2}b_{2}\right) ^{2}\left( b_{1}c_{11}+b_{2}c_{21}\right)
+\left( a_{1}\right) ^{2}b_{1}b_{2}\left( b_{1}c_{14}+b_{2}c_{24}\right)
\nonumber \\
&&\left. +a_{1}a_{2}b_{1}\left[ \left( b_{1}\right) ^{2}c_{14}-\left(
b_{2}\right) ^{2}c_{23}-b_{1}b_{2}\left( c_{13}-c_{24}\right) \right]
\right\} /\left[ a_{1}\left( a_{2}\right) ^{2}c^{2}\right] ~.
\label{333gam3}
\end{eqnarray}

The fourth choice are the $3$ eqs. (\ref{33K2}), (\ref{33K3}) and (\ref{33K4}%
). One then gets (if need be, with the help of \textbf{Mathematica)}
\end{subequations}
\begin{subequations}
\label{334gamj}
\begin{eqnarray}
&&\gamma _{1}=\left\{ -3\left( a_{2}b_{1}\right)
^{2}b_{2}+6a_{1}a_{2}b_{1}(b_{2})^{2}-3\left( a_{1}\right) ^{2}\left(
b_{2}\right) ^{3}\right.  \nonumber \\
&&\left. +\left( a_{2}\right) ^{2}\left( b_{1}c_{12}+b_{2}c_{22}\right)
-2a_{1}a_{2}\left( b_{1}c_{13}+b_{2}c_{23}\right) \right.  \nonumber \\
&&\left. +3\left( a_{1}\right) ^{2}\left( b_{1}c_{14}+b_{2}c_{24}\right)
\right\} /\left( a_{2}c^{2}\right) ~,  \label{334gam1}
\end{eqnarray}%
\begin{eqnarray}
&&\gamma _{2}=\left\{ 3\left( a_{2}b_{1}b_{2}\right) ^{2}+3\left(
a_{1}\right) ^{2}\left( b_{2}\right) ^{4}-6a_{1}a_{2}b_{1}\left(
b_{2}\right) ^{3}\right.  \nonumber \\
&&\left. -3\left( a_{1}\right) ^{2}b_{2}\left(
b_{1}c_{14}+b_{2}c_{24}\right) \right.  \nonumber \\
&&\left. +\left( a_{2}\right) ^{2}\left[ \left( b_{1}\right)
^{2}c_{13}-2\left( b_{2}\right) ^{2}c_{22}+b_{1}b_{2}\left(
-2c_{12}+c_{23}\right) \right] \right.  \nonumber \\
&&\left. -3a_{1}a_{2}\left[ \left( b_{1}\right) ^{2}c_{14}-\left(
b_{2}\right) ^{2}c_{23}+b_{1}b_{2}\left( -c_{13}+c_{24}\right) \right]
\right\} /\left( a_{2}c\right) ^{2}~,  \label{334gam2}
\end{eqnarray}

\begin{eqnarray}
&&\gamma _{3}=\left\{ -\left( a_{1}\right) ^{2}\left( b_{2}\right)
^{5}-\left( a_{2}\right) ^{2}\left( b_{1}\right) ^{2}\left( b_{2}\right)
^{3}+2a_{1}a_{2}b_{1}\left( b_{2}\right) ^{4}\right.  \nonumber \\
&&+\left( a_{1}b_{2}\right) ^{2}\left( b_{1}c_{14}+b_{2}c_{24}\right)
+\left( a_{2}\right) ^{2}\left[ \left( b_{2}\right) ^{3}c_{22}+\left(
b_{1}\right) ^{3}c_{14}\right]  \nonumber \\
&&+\left( a_{2}\right) ^{2}b_{1}b_{2}\left[ b_{1}\left(
-c_{13}+c_{24}\right) +b_{2}\left( c_{12}-c_{23}\right) \right]  \nonumber \\
&&\left. +a_{1}a_{2}b_{2}\left[ \left( b_{1}\right) ^{2}c_{14}-\left(
b_{2}\right) ^{2}c_{23}+b_{1}b_{2}\left( -c_{13}+c_{24}\right) \right]
\right\} /\left[ \left( a_{2}\right) ^{3}c^{2}\right] ~.  \label{334gam3}
\end{eqnarray}

The formulas displayed above provide \textit{explicit} expressions of the $3$
parameters $\gamma _{j}$ in terms of the $8$ coefficients $c_{n\ell },$ and
also in terms of the parameters $a_{n}$ and $b_{n}$ the determination of
which in terms of the $8$ coefficients $c_{n\ell }$ has been detailed in the
preceding subsections of this \textbf{Section 3}.

These $4$ expressions of the $3$ parameters $\gamma _{j}$ are of course
equivalent provided the $8$ coefficients $c_{n\ell }$ satisfy the $2$
\textit{constraints} determined above (see \textbf{Subsection 3.2} and
\textbf{Appendices A}, \textbf{B} and \textbf{C}).

The task of inverting the formulas (\ref{2}) is thereby completed.

\subsection{The special case with $c_{14}=c_{21}=0$}

In several applications it is \textit{unreasonable} to assume that the
variations over time of the quantity $x_{n}\left( t\right) $ be influenced
by a cause which depends \textit{only} on the value of the \textit{other}
variable:\textit{\ }hence that in the right-hand sides of the ODE (\ref{1})
characterizing the change over time of the dependent variable $x_{n}\left(
t\right) $ associated to these phenomena, terms \textit{independent} from
the values of this variable $x_{n}\left( t\right) $ be present. This fact
motivates the special interest of the \textit{subclass} of the dynamical
systems (\ref{1}) characterized by the vanishing of the $2$ coefficients $%
c_{14}$ and $c_{21},$%
\end{subequations}
\begin{subequations}
\begin{equation}
c_{14}=c_{21}=0~,  \label{35Con}
\end{equation}%
hence characterized by the following reduced version of the system (\ref{1}):%
\begin{eqnarray}
&&\dot{x}_{1}=x_{1}\left[ g_{11}\left( x_{1}\right)
^{2}+g_{12}x_{1}x_{2}+g_{13}\left( x_{2}\right) ^{2}\right] ~,  \nonumber \\
&&\dot{x}_{2}=x_{2}\left[ g_{21}\left( x_{1}\right)
^{2}+g_{22}x_{1}x_{2}+g_{23}\left( x_{2}\right) ^{2}\right] ~,
\label{5x12dot}
\end{eqnarray}%
where, for notational convenience, we set (in addition to (\ref{35Con}))%
\begin{equation}
c_{1j}\equiv g_{1j}~,~~~c_{2\ell }\equiv g_{2j}~,~~~j=1,2,3~,~~~\ell =j+1~.
\label{gggg}
\end{equation}%
In this \textbf{Subsection 3.5 }we tersely treat this particular \textit{%
subcase} which---for the reason mentioned above---is of special \textit{%
applicative} relevance.

\textbf{Remark 3.5-1}. Clearly in this case the invariance property (\ref%
{1Change}) is replaced by the following formulas:
\end{subequations}
\begin{equation}
x_{1}\left( t\right) \Leftrightarrow x_{2}\left( t\right)
~;~g_{11}\Leftrightarrow g_{23}~,~g_{12}\Leftrightarrow
g_{22}~,~g_{13}\Leftrightarrow g_{21}~.  \label{35Trans}
\end{equation}%
While of course for the system (\ref{5x12dot}) the invariance property (\ref%
{1Inveta}) is just as valid as for the system (\ref{1}). $\blacksquare $

The most direct way to obtain the results relevant to this special case is
to insert in the previous treatment the restriction (\ref{35Con}) and the
notational change (\ref{gggg}).

To assess the impact of the restriction (\ref{35Con}) we insert
it---together with the notational change (\ref{gggg})---in the $2$ eqs. (\ref%
{c14c21}).

It is then easily seen that its insertion in the simpler eqs. (\ref{c14c21a}%
) implies the restriction
\begin{subequations}
\begin{equation}
g_{12}=g_{22}=0~,
\end{equation}%
causing the system (\ref{5x12dot}) to take the following reduced form:%
\begin{eqnarray}
&&\dot{x}_{1}=x_{1}\left[ g_{11}\left( x_{1}\right) ^{2}+g_{13}\left(
x_{2}\right) ^{2}\right] ~,  \nonumber \\
&&\dot{x}_{2}=x_{2}\left[ g_{21}\left( x_{1}\right) ^{2}+g_{23}\left(
x_{2}\right) ^{2}\right] ~;
\end{eqnarray}%
and moreover implying---see (\ref{32a2})---the vanishing of the parameter $%
a_{2},$%
\begin{equation}
a_{2}=0~,
\end{equation}
signifying that in this special case our solution technique gets applied to
a system whose solvable character is rather trivial.

The situation is instead different if one considers the alternative case (%
\ref{c14c21b}). Then, via (\ref{35Con}) and (\ref{gggg}), one obtains the
following $2$ \textit{new constraints} on the $6$ coefficients $g_{nj}$:
\end{subequations}
\begin{subequations}
\label{35NewCon}
\begin{eqnarray}
-g_{12}\left[ 2\left( g_{22}\right) ^{2}+g_{21}\left( g_{13}-3g_{23}\right) %
\right] &&  \nonumber \\
+g_{22}\left[ \left( g_{12}\right) ^{2}+\left( g_{22}\right)
^{2}+3g_{21}\left( g_{13}-g_{23}\right) \right] +FR_{4}=0~, &&
\end{eqnarray}%
\begin{eqnarray}
-\left( g_{12}\right) ^{4}g_{22}+g_{13}\left( 3g_{11}-g_{21}\right) \left\{
g_{22}\left[ \left( g_{22}\right) ^{2}+g_{21}\left( g_{13}-3g_{23}\right) %
\right] -FR_{4}\right\} &&  \nonumber \\
+\left( g_{12}\right) ^{3}\left[ 3\left( g_{22}\right) ^{2}+g_{21}\left(
g_{13}-3g_{23}\right) \right] +g_{12}\left\{ -3\left( g_{13}\right)
^{2}g_{21}\left( g_{11}-g_{21}\right) +\left( g_{22}\right) ^{4}\right. &&
\nonumber \\
\left. -g_{22}FR_{4}-3g_{21}\left( g_{22}\right) ^{2}g_{23}+g_{13}\left[
9g_{21}g_{23}\left( g_{11}-g_{21}\right) -\left( g_{22}\right) ^{2}\left(
6g_{11}-5g_{21}\right) \right] \right\} &&  \nonumber \\
+\left( g_{12}\right) ^{2}\left\{ -3\left( g_{22}\right) ^{3}+FR_{4}+g_{22}
\left[ 3g_{11}g_{13}-g_{21}\left( 5g_{13}-6g_{23}\right) \right] \right\}
=0~, &&
\end{eqnarray}%
with
\begin{equation}
FR_{4}=\left[ g_{21}\left( g_{13}-3g_{23}\right) -g_{22}\left(
g_{12}-g_{22}\right) \right] \sqrt{\left( g_{12}-g_{22}\right)
^{2}+4g_{13}g_{21}}~.
\end{equation}

These $2$ \textit{constraints} on the $6$ coefficients $g_{nj}$ are of
course required to hold---for the applicability of the solution of the
system (\ref{5x12dot}) via our approach---in addition to the other $2$
\textit{constraints} on the $6$ coefficients $g_{nj}$ inherited from the
previous treatment, which may be obtained via (\ref{35Con}) and (\ref{gggg})
from any pair of the formulas of \textbf{Appendix C}.

An alternative---and perhaps more convenient---way to obtain alternative,
perhaps more useful, versions of these \textit{constraints} is by revisiting
the treatment of the preceding \textbf{Subsections} of this \textbf{Section 3%
}, highlighting the modifications implied by the restriction (\ref{35Con})
together with the notational change (\ref{gggg}).

The $3$ expressions (\ref{32alfa}) of the parameter $\alpha $ read then as
follows:
\end{subequations}
\begin{subequations}
\label{5alpha3}
\begin{equation}
\alpha =\frac{g_{12}}{3g_{11}-g_{21}}~,  \label{5alpha}
\end{equation}%
~%
\begin{equation}
\alpha =\frac{g_{13}g_{22}}{g_{12}g_{22}-g_{13}g_{21}-\left( g_{22}\right)
^{2}+3g_{21}g_{23}}~,
\end{equation}

\begin{equation}
\alpha =\frac{g_{12}}{3g_{11}-g_{21}}~;
\end{equation}%
and (since the first and last of these $3$ equations are clearly identical)
they clearly yield the following \textit{single constraint} on the $6$
coefficients $g_{nj}$ ($n=1,2$; $j=1,2,3$):
\end{subequations}
\begin{equation}
g_{12}\left[ g_{12}g_{22}-g_{13}g_{21}-\left( g_{22}\right)
^{2}+3g_{21}g_{23}\right] -g_{13}g_{22}\left( 3g_{11}-g_{21}\right) =0~.
\label{35Constraint}
\end{equation}

The nonlinear context implies that there are \textit{many equivalent}
versions of this \textit{constraint}, and as well \textit{many equivalent}
versions of a \textit{second} independent relation \textit{constraining} the
coefficients $g_{nj},$ the simultaneous validity of both of them being then
required in order that the solution of the system (\ref{5x12dot}) be
provided by the insertion in \textbf{Proposition 1-1} of the identities (\ref%
{35Con}) and (\ref{gggg}). The identification of which version of these $2$
\textit{constraints} to be used is not obvious. Perhaps the most convenient
versions are provided by the following formulas:
\begin{subequations}
\label{DCons}
\begin{equation}
-3g_{11}g_{13}g_{22}+\left( g_{12}-g_{22}\right) \left(
g_{12}g_{22}-g_{13}g_{21}\right) +3g_{12}g_{21}g_{23}=0~,
\end{equation}%
\begin{equation}
g_{12}g_{22}\left( g_{13}-3g_{23}\right) \left( g_{21}-3g_{11}\right)
-\left( g_{12}g_{22}\right) ^{2}=0~,
\end{equation}%
obtained from the relations (\ref{A1}) and (\ref{32SecCona}) via the
relations (\ref{35Con}) and (\ref{gggg}).

Next---for the convenience of the interested reader---let us report the
expressions of \textit{every pair} of the $6$ coefficients $g_{nj}$ in terms
of the other $4$ implied by these $2$ constraints.

There are of course altogether $15$ different pairs of the $6$ coefficients $%
g_{nj}$ ($15=6\cdot 5/2$), but only $9$ relations expressing each pair are
reported below, since the additional $6$ relations are implied by the last $%
6 $ reported below via the transformations (\ref{35Trans}). All these
expressions can be obtained by solving simultaneously the $2$ \textit{%
constraints} (\ref{DCons}); remarkably, this task can be implemented \textit{%
explicitly} (with the help of \textbf{Mathematica}) and the outcomes are
quite simple.

\textbf{Remark 3.5-2}. The formulas reported below are valid for \textit{%
generic} values of the coefficients $g_{nj}$: alternative solutions
requiring some of the coefficients $g_{nj}$ to \textit{vanish} are generally
\textit{not} reported, unless they are exceptionally simple (moreover, in
the last case reported below, see (\ref{35ggf}), it is the \textit{only}
result obtained for that case). $\blacksquare $

Hereafter
\end{subequations}
\begin{eqnarray}
G_{1} &=&\sqrt{\left( g_{12}\right) ^{2}+4g_{13}g_{21}-2g_{12}g_{22}+\left(
g_{22}\right) ^{2}}~,  \nonumber \\
G_{11/23} &=&\sqrt{g_{11}/g_{23}}~,  \nonumber \\
G_{2} &=&\sqrt{\left( g_{12}\right) ^{2}+12g_{21}g_{23}}~.  \label{GGG}
\end{eqnarray}

We first report the formulas expressing the $3$ pairs that are transformed
into themselves by the transformations (\ref{35Trans}):
\begin{subequations}
\begin{equation}
g_{11}=\left[ \left( g_{12}\right)
^{2}+2g_{13}g_{21}-g_{12}g_{22}+g_{12}G_{1}\right] /\left( 6g_{13}\right) ~,
\end{equation}%
\begin{equation}
g_{23}=\left[ 2g_{13}g_{21}-g_{12}g_{22}+\left( g_{22}\right)
^{2}+g_{12}G_{1}\right] /\left( 6g_{21}\right) ~;
\end{equation}%
\begin{equation}
g_{12}=-\left( 3g_{11}-g_{21}\right) /G_{11/23}~,~~~g_{22}=\left(
g_{13}-3g_{23}\right) G_{11/23}~~;
\end{equation}%
\begin{equation}
g_{13}=3g_{23}-g_{22}/G_{11/23}~,~~~g_{21}=3g_{11}-g_{12}G_{11/23}~.
\end{equation}

Next, the formulas expressing the $6$ pairs that are transformed into $6$
\textit{other} pairs by the transformations (\ref{35Trans}):
\end{subequations}
\begin{subequations}
\label{35ggabcdef}
\begin{eqnarray}
g_{11} &=&g_{21}/3~,~~~g_{12}=0~;  \nonumber \\
g_{11} &=&\frac{\left( g_{22}\right) ^{2}g_{23}}{\left(
g_{13}-3g_{23}\right) ^{2}}~,  \nonumber \\
g_{12} &=&\frac{\left( g_{13}\right) ^{2}g_{21}-6g_{13}g_{21}g_{23}-3\left[
\left( g_{22}\right) ^{2}g_{23}-3g_{21}\left( g_{23}\right) ^{2}\right] }{%
g_{22}\left( g_{13}-3g_{23}\right) }~;  \label{35gga}
\end{eqnarray}%
\begin{equation}
g_{11}=\frac{\left( g_{12}\right) ^{2}+6g_{21}g_{23}+g_{12}G_{2}}{18g_{23}}%
~,~~~g_{13}=\frac{g_{12}g_{22}+6g_{21}g_{23}-g_{22}G_{2}}{2g_{21}}~;
\label{35ggb}
\end{equation}%
\begin{equation}
g_{11}=\frac{\left( g_{22}\right) ^{2}g_{23}}{\left( g_{13}-3g_{23}\right)
^{2}}~,~~~g_{21}=\frac{g_{22}\left[ g_{12}\left( g_{13}-3g_{23}\right)
+3g_{22}g_{23}\right] }{\left( g_{13}-3g_{23}\right) ^{2}}~;  \label{35ggc}
\end{equation}

\begin{equation}
g_{11}=\frac{6g_{21}g_{23}+g_{12}\left( g_{12}+G_{2}\right) }{18g_{23}}%
~,~~~g_{22}=\frac{-\left( g_{13}-3g_{23}\right) \left( g_{12}+G_{2}\right) }{%
6g_{23}}~;  \label{35ggd}
\end{equation}%
\begin{equation}
g_{12}=\left( 3g_{11}-g_{21}\right)
/G_{11/23}~,~~~g_{13}=3g_{23}+g_{22}/G_{11/23}~;  \label{35gge}
\end{equation}%
\begin{equation}
g_{12}=0~,~~~g_{21}=3g_{11}~.  \label{35ggf}
\end{equation}

Hereafter we assume of course that the $6$ coefficients $g_{nj}$ satisfy the
$2$ \textit{constraints}, many versions of which are provided above; and in
addition the $2$ constraints (\ref{35NewCon}). Hence all the formulas we
display below (in this \textbf{Subsection 3.5}) are just representative
\textit{avatars} of many other \textit{equivalent} expressions implied by
the formulas reported above, relating the $6$ coefficients $g_{nj}$ to each
others.

The solution of the initial-values problem for the system (\ref{5x12dot}) is
then provided by \textbf{Proposition 1-1}, complemented by the eqs. (\ref%
{35Con}) and (\ref{gggg}).

For the applications of these findings the "inverse problem" to express the $%
7$ parameters $a_{n},$ $b_{n}$ and $\gamma _{j}$ in terms of the
coefficients $g_{nj}$ is of course also quite important.

A quite neat version of the formulas expressing the $2$ parameters $a_{n}$
in terms of the coefficients $g_{nj}$ of the system (\ref{5x12dot}) reads as
follows:
\end{subequations}
\begin{equation}
a_{1}=\sqrt{g_{11}}~,~~~a_{2}=\sqrt{g_{23}}~.
\end{equation}

For the determination of the $2$ parameters $b_{n}$ we refer to the relevant
treatment provided at the end of \textbf{Subsection 3.3}, which is
applicable with the following modifications: the equations (\ref{betaeta})
read now
\begin{equation}
a_{2}=\alpha a_{1}~,~~~g_{nj}=\eta \hat{g}_{nj}~,~~~b_{1}=\eta
~,~~~b_{2}=\alpha \beta \eta ~,  \label{35betaeta}
\end{equation}%
where now $\alpha $ is given in terms of the coefficient $g_{nj}$ by anyone
of the formulas (\ref{5alpha3}); $\eta $ is a parameter characterizing now a
rescaling of the $6$ coefficients $g_{nj}$---introduced here to make
notational contact with the invariance property (see (\ref{1Inveta})) of the
system (\ref{5x12dot}) mentioned above (see \textbf{Remark 3.5-1 }and\textbf{%
\ Remark 1-1}), again implying that this parameter can be eventually
altogether eliminated---i. e., assigned an \textit{arbitrary} \textit{%
nonvanishing value }(for instance, just the value $\eta =1$)---via a
corresponding appropriate rescaling of the independent variable $t$; $\beta $
is the parameter that we determine immediately below in terms of the $6$
coefficients $g_{nj}$; and the last $2$ eqs. (\ref{35betaeta}) determine of
course the $2$ parameters $b_{1}$ and $b_{2}$, thereby completing the task
indicated by the title of this \textbf{Subsection 3.3}.

To determine the parameter $\beta $ in terms of the $6$ parameters $g_{nj}$%
---or, equivalently, $\hat{g}_{n\ell }$ (see (\ref{35betaeta}) and (\ref%
{1Inveta}))---we insert in eq. (\ref{33bb}) the positions (\ref{35betaeta}),
of course with (\ref{35Con}) and (\ref{gggg}), getting thereby the following
\textit{cubic} equation for this parameter:
\begin{equation}
\left( \beta -1\right) ^{3}=\beta \left( \hat{g}_{23}\alpha ^{-2}-\hat{g}%
_{22}\alpha ^{-1}+\hat{g}_{21}\right) -\hat{g}_{13}\alpha ^{-1}+\hat{g}_{12}-%
\hat{g}_{11}\alpha ^{-1}~.  \label{35beta}
\end{equation}%
\textit{Explicit} solutions of this equation can of course be provided via
the standard Cardano formulas.

Finally, several equivalent formulas expressing the parameters $\gamma _{j}$
and involving the parameters $a_{n}$ and $b_{n}$ defined just above are
provided by the formulas (\ref{331gamj})-(\ref{334gamj}), of course after
the replacement implied by the identities (\ref{35Con}) and (\ref{gggg}) of
the $8$ coefficients $c_{n\ell }$ with the $6$ coefficients $g_{nj}$ (or
possibly $\hat{g}_{nj}$, see (\ref{35betaeta})).

\section{The isochronous extension}

In this \textbf{Section 4} attention is restricted to \textit{real} values
of the independent variable $t$ ("time"); except for the discussion in
\textbf{Remark 4.3}.

It is well known (see, if need be, \cite{FC2008} and references therein)
that via the following simple change of dependent and independent variables,
\begin{subequations}
\begin{equation}
\tilde{x}_{n}\left( t\right) =\exp \left( \mathbf{i}\omega t\right)
x_{n}\left( \tau \right) ~,~~~\tau =\left[ \exp \left( 2\mathbf{i}\omega
t\right) -1\right] /\left( 2\mathbf{i}\omega \right) ~,~~~n=1,2~,
\label{VarChange}
\end{equation}%
the (\textit{autonomous}) system (\ref{1}) gets transformed into the, also
\textit{autonomous}, system (\ref{1Iso}). This new system differs from the
system (\ref{1}) only due to the additional presence of the linear term $%
\mathbf{i}\omega \tilde{x}_{n}$ in the right-hand side of its $2$ ODEs. It
is of course just as \textit{solvable} as the system (\ref{1}), to which it
gets reduced via the (easily invertible) change of variables (\ref{VarChange}%
). One then notes (see, for instance, \cite{FC2008}) that, if the parameter $%
\mathbf{i}\omega $ is an \textit{arbitrary imaginary} number (as we
hereafter assume), then---see (\ref{VarChange})---as the (\textit{real})
time variable $t$ evolves from its initial value $0$ towards $+\infty ,$ the
auxiliary \textit{complex} variable $\tau $ rotates (counterclockwise for $%
\omega >0,$ clockwise for $\omega <0$) on the circle $C$ of radius $%
1/\left\vert 2\omega \right\vert $ centered at the point $\mathbf{i/}\left(
2\omega \right) $ in the \textit{complex} $\tau $-plane: a clearly \textit{%
periodic} evolution, with period%
\begin{equation}
T=\pi /\left\vert \omega \right\vert ~.  \label{4T}
\end{equation}%
Hence any analytic function $f\left( \tau \right) $ of the \textit{complex}
variable $\tau $ featuring no singularity \textit{inside} (nor \textit{on})
the circle $C$ in the complex $\tau $-plane evolves \textit{periodically}
with that period $T$ as function of \textit{time}, namely as the function $%
\tilde{f}\left( t\right) \equiv f\left( \tau \left( t\right) \right) $ of
the \textit{real} variable $t$ ("time"); and it also evolves \textit{%
periodically in time}---with a period which is then an \textit{integer
multiple} of $T$---if the function $f\left( \tau \right) $ features, as
\textit{analytic} function of the \textit{complex} variable $\tau $, a
\textit{finite} number of \textit{rational} \textit{branch points} inside
the circle $C$ (except for \textit{nongeneric initial data} such that one of
these branch points falls exactly on the circle $C$ in the \textit{complex} $%
\tau $-plane); since then the argument $\tau \left( t\right) $ of the
function $\tilde{f}\left( t\right) $ travels on a Riemann surface with a
\textit{finite} number of sheets, and while possibly visiting (some of)
these sheets it always eventually retraces the same path on the Riemann
surface. Hence the system (\ref{1Iso})---considered as a function of \textit{%
real} time $t$---has the remarkable feature to be \textit{isochronous}
whenever the corresponding system (\ref{1}) belongs to the class identified
in this paper, such that its solutions only feature a \textit{finite} number
of \textit{rational} branch points (see \textbf{Remark 1-4}): the period of
its solutions being then generally a \textit{finite integer multiple} of the
basic period $T,$ see (\ref{4T}), which does not change for sufficiently
small changes of the initial data $\tilde{x}_{n}\left( 0\right) =x_{n}\left(
0\right) $ (see (\ref{VarChange})); it may discontinuously change when the
initial data are modified so that the consequential shifts of the branch
point positions in the complex $\tau $-plane causes one of them to enter or
exit the circle $C$. While the special solutions---characterized by the
\textit{exceptional initial data }$\tilde{x}_{1}\left( 0\right) $ and $%
\tilde{x}_{2}\left( 0\right) $ such that one of the branch points of the
corresponding solutions $\tilde{x}_{n}\left( t\right) $ falls \textit{exactly%
} on the circle $C$---then feature---at some special value (or values) $%
t_{s} $ of the time variable $t$---a coincidence of the $2$ functions $%
\tilde{x}_{1}\left( t\right) $ and $\tilde{x}_{2}\left( t\right) $, $\tilde{x%
}_{1}\left( t_{s}\right) =\tilde{x}_{2}\left( t_{s}\right) ,$ namely a
"collision" of the $2$ \textit{complex} points $\tilde{x}_{1}\left( t\right)
$ and $\tilde{x}_{2}\left( t\right) $, causing thereafter a loss of their
individual identities, a phenomenon whose occurrence is associated with the
singularity of that particular solution of the system (\ref{1Iso}). But of
course these phenomena only happen for \textit{nongeneric} values of the
initial data $\tilde{x}_{1}\left( 0\right) $ and $\tilde{x}_{2}\left(
0\right) $.

\textbf{Remark 4-1}. Actually the periods of the \textit{generic} solutions $%
\tilde{x}_{n}\left( t\right) $ of the system (\ref{1Iso})---whenever the
corresponding solutions of the system (\ref{1}) are \textit{algebraic}---are
\textit{integer multiples} of the period $\tilde{T}=2T=2\pi /\left\vert
\omega \right\vert $, due to the prefactor $\exp \left( \mathbf{i}\omega
t\right) =\exp \left( 2\mathbf{i}\pi st/\tilde{T}\right) $ (where $s=\omega
/\left\vert \omega \right\vert $ is the sign of $\omega $) in the change of
variables from $x_{n}\left( \tau \right) $ to $\tilde{x}_{n}\left( t\right)
, $ see (\ref{VarChange}) and (\ref{4T}); while $x_{n}\left( \tau \right)
\equiv x_{n}\left( \tau \left( t\right) \right) $ is clearly itself periodic
in $t$ with period $T$ (hence as well with period $\tilde{T}=2T$) if the
functions $x_{n}\left( \tau \right) $ are \textit{holomorphic} in $\tau $,
and instead with a \textit{positive integer} multiple of $T$ if the
functions $x_{n}\left( \tau \right) $ are \textit{not holomorphic} but only
feature a \textit{finite} number of \textit{rational} branch points as
functions of the \textit{complex} variable $\tau $. $\blacksquare $

The possibility reviewed in this \textbf{Section 4} (and see \cite{FC2008}
for analogous treatments in more general contexts), to associate to a
\textit{homogeneous} dynamical system such as (\ref{1}) a corresponding
\textit{isochronous} system such as (\ref{1Iso}), underlines the interest to
identify classes of \textit{homogeneous} dynamical systems the general
solutions of which---when considered as \textit{analytic} functions of
\textit{complex} time---only feature a \textit{finite} number of \textit{%
rational} branch points; and in some cases to even obtain \textit{explicitly}
the solutions of their initial-values problems.

\textbf{Remark 4-2}. Let us also emphasize the obvious fact that the \textit{%
isochronous} systems obtained in this manner---such as (\ref{1Iso})---entail
a \textit{doubling} of the (\textit{real}) dependent variables, since the
presence of the \textit{imaginary} parameter $\mathbf{i}\omega $ in their
right-hand side entails the need to treat the dependent variables $\tilde{x}%
_{n}\left( t\right) $ as \textit{complex} numbers, featuring both a \textit{%
real} and an \textit{imaginary} part,
\end{subequations}
\begin{equation}
\tilde{x}_{n}\left( t\right) \equiv \tilde{x}_{Rn}\left( t\right) +\mathbf{i}%
\tilde{x}_{In}\left( t\right) ~,~~~n=1,2~;
\end{equation}%
hence the $2$ ODEs (\ref{1Iso}) define a dynamical system involving \textit{%
de facto} the $4$ \textit{real} dependent variables $\tilde{x}_{R1}\left(
t\right) ,$ $\tilde{x}_{R2}\left( t\right) ,$ $\tilde{x}_{I1}\left( t\right)
,$ $\tilde{x}_{I2}\left( t\right) ;$ and of course an analogous \textit{%
doubling by complexification} may conveniently be made for its coefficients $%
c_{n\ell }\equiv c_{Rn\ell }+\mathbf{i}c_{In\ell }$ whenever planning to use
it in an \textit{applicative} context, which generally entails the use of
\textit{real numbers}. $\blacksquare $

Let us complete this \textbf{Section 4 }by reporting (see for instance \cite%
{FC2008} and references therein) the following well-known

\textbf{Remark 4-3}. It is plain that, for\textit{\ any arbitrary}
assignment of the $8$ parameters $c_{n\ell },$ the data $x_{1}=x_{2}=0$
identify an \textit{equilibrium} configuration of the system (\ref{1}); and
moreover that if the initial data $x_{1}\left( 0\right) $ and $x_{2}\left(
0\right) $ are \textit{sufficiently small}---say, $\left\vert x_{1}\left(
0\right) \right\vert <\varepsilon $ and $\left\vert x_{2}\left( 0\right)
\right\vert <\varepsilon $ with $\varepsilon $ a sufficiently small
parameter---then the corresponding solution $x_{n}\left( t\right) $ of the
system (\ref{1}) features the property to be \textit{holomorphic} as a
function of its (\textit{complex}) argument $t$ in the neighborhood of the
initial data---namely, there exists a parameter $\delta $ (possibly quite
small, but \textit{strictly positive}, $\delta >0$) such that $x_{n}\left(
t\right) $ are both \textit{holomorphic} functions of the \textit{complex}
variable $t$ inside the disk $\left\vert t\right\vert <\delta .$ This
clearly implies that the general system (\ref{1Iso})---for \textit{any
arbitrary} assignment of its $8$ parameters $c_{n\ell }$ and a sufficiently
\textit{small} value of the (\textit{real, nonvanishing}) parameter $\omega $%
---features the \textit{isochronicity} property to be---as function of the
\textit{real} variable $t$ ("time")---\textit{completely periodic} with
period $\tilde{T}=2T$ (see (\ref{4T})) for the \textit{open} set of its
initial data $\tilde{x}_{1}\left( 0\right) =x_{1}\left( 0\right)
<\varepsilon $ and $\tilde{x}_{2}\left( 0\right) =x_{2}\left( 0\right)
<\varepsilon $ such that $2\left\vert \omega \right\vert \delta <1$. $%
\blacksquare $

And let us also mention that the \textit{isochronous extension} described in
this \textbf{Section 4} is of course also applicable to the special case
treated in \textbf{Subsection 3.5}; the relevant treatment is left as an
easy exercise for the interested reader.

\section{Conclusions and outlook}

In this last section we tersely review the main findings of this paper, and
then mention possible analogous developments.

The \textit{initial-values} problem of the system (\ref{1}) can be
solved---for \textit{arbitrary} initial data $x_{1}\left( 0\right) $ and $%
x_{2}\left( 0\right) $---provided the $8$ \textit{a priori arbitrary}
coefficients $c_{n\ell }$ satisfy $2$ \textit{constraints}: as discussed in
detail above, see in particular the treatment in \textbf{Section 3} and the
determination of $2$ of the $8$ coefficients $c_{n\ell }$ in terms of the
other $6$ implied by these $2$ \textit{constraints}, as reported in \textbf{%
Appendix C}. The solution is provided by \textbf{Proposition 1-1}. It is
given \textit{rather explicitly} there in terms of the $7$ parameters $%
a_{1}, $ $a_{2},$ $b_{1},$ $b_{2},$ $\gamma _{1},$ $\gamma _{2},$ $\gamma
_{3},$ themselves given in terms of the $8$ coefficients $c_{n\ell },$ as
follows: the $2$ parameters $a_{1}$ and $a_{2}$ are \textit{explicitly}
given, for instance, by eqs. (\ref{32aa12}) with $\alpha $ given by any one
of the $3$ formulas (\ref{32alfa}); the $2$ parameters $b_{1},$ $b_{2}$ can
be determined as explained at the end of \textbf{Subsection 3.3}; and the $3$
parameters $\gamma _{1},$ $\gamma _{2},$ $\gamma _{3}$ are then defined by
any one of the $4$ sets of \textit{explicit} formulas displayed in the last
part of \textbf{Subsection 3.4}. The remaining ambiguities in these
determinations are then eventually eliminated when the solution of the
\textit{initial-values} problem is expressed---see \textbf{Proposition 1-1}%
---in terms of the $2$ initial data $x_{1}\left( 0\right) $ and $x_{2}\left(
0\right) $; and the solutions $x_{1}\left( t\right) $ and $x_{2}\left(
t\right) $ are then identified by \textit{continuity} in $t$ starting from
the initial data $x_{1}\left( 0\right) $ and $x_{2}\left( 0\right) $.

In an analogous manner the \textit{initial-values} problem can be solved for
the related \textit{isochronous} system of ODEs (\ref{1Iso}); while the
conditions implying its \textit{isochrony} are detailed in \textbf{Remark 1-4%
} complemented by the relevant formulas of \textbf{Proposition 1-1} and by
the treatment of this case in \textbf{Section 4 }(including in particular
\textbf{Remarks 4-2} and \textbf{4-3)}.

Finally, let us mention that natural extensions of the findings reported in
this paper can be sought by considering dynamical systems analogous to, but
more general than, (\ref{1}): in particular systems involving more than $2$
dependent variables, and systems featuring in their right-hand sides \textit{%
higher-degree}---and possibly \textit{nonhomogeneous}---polynomials; as well
as dynamical systems involving \textit{higher-order} ODEs; or evolutions in
\textit{discrete}---rather than \textit{continuous}---time.

\section{Acknowledgements}

It is a pleasure to thank Fran\c{c}ois Leyvraz and Robert Conte for useful
suggestions. One of us (FP) would like to thank the Physics Department of
the University of Rome "La Sapienza" for hosting her in the Department
during several visits to Rome, and Payame Noor University (PNU) for
financial support to those visits during her sabbatical and for a grant for
this research; and the Italian Embassy in Iran as well as the Italian
Ministry of Foreign Affairs for two grants for short visits to Italy
("contributi finalizzati a brevi visite di docenti e di personalita' della
cultura straniere in Italia").

\bigskip

\section{Appendix A}

In this \textbf{Appendix A} the formulas ---implied by the \textit{first
constraint} on the $8$ coefficients $c_{n\ell }$ obtained in \textbf{%
Subsection 3.2}: see for instance eq. (\ref{32Constraint2}) or any one of
the $3$ equivalent eqs. (\ref{32Constraints})---are reported, which express
each of these $8$ coefficients in terms of the other $7$; only $4$ of these
expressions are actually displayed below, since the other $4$ can be
obtained from those displayed via the transformations (\ref{1Change}). These
formulas provide of course $8$ \textit{constraints} on the $8$ coefficients $%
c_{n\ell }$. These \textit{constraints} are \textit{all} \textit{equivalent
among themselves}.
\begin{subequations}
\label{AAA}
\begin{eqnarray}
c_{11} &=&\left\{ -3\left( c_{13}\right) ^{2}c_{21}-\left( c_{12}\right)
^{2}c_{23}+3c_{14}\left[ \left( c_{22}\right) ^{2}-3c_{21}c_{23}\right]
\right.  \nonumber \\
&&\left. +c_{13}\left[ c_{22}(c_{12}-c_{23})+9c_{21}c_{24}\right] +c_{12}%
\left[ 9c_{14}c_{21}+\left( c_{23}\right) ^{2}-3c_{22}c_{24}\right] \right\}
/  \nonumber \\
&&/\left[ 3\left( 3c_{14}c_{22}-c_{13}c_{23}\right) \right] ~,  \label{A1}
\end{eqnarray}%
\begin{eqnarray}
&&c_{12}=\left[ 9c_{14}c_{21}+c_{13}c_{22}+\left( c_{23}\right)
^{2}-3c_{22}c_{24}\right.  \nonumber \\
&&\pm \left\{ \left[ 9c_{14}c_{21}+c_{13}c_{22}+\left( c_{23}\right)
^{2}-3c_{22}c_{24}\right] ^{2}\right.  \nonumber \\
&&-4c_{23}\left[ 3\left( c_{13}\right)
^{2}c_{21}+9c_{11}c_{14}c_{22}-3c_{14}\left( c_{22}\right)
^{2}+9c_{14}c_{21}c_{23}\right.  \nonumber \\
&&\left. \left. \left. +c_{13}\left(
-3c_{11}c_{23}+c_{22}c_{23}-9c_{21}c_{24}\right) \right] \right\} ^{1/2}%
\right] /\left( 2c_{23}\right) ~,  \label{A2}
\end{eqnarray}%
\begin{eqnarray}
c_{13} &=&\left[ 9c_{21}c_{24}+c_{12}c_{22}+3c_{11}c_{23}-c_{22}c_{23}\right.
\nonumber \\
&&\pm \left\{ \left[ \left(
9c_{21}c_{24}+c_{12}c_{22}+3c_{11}c_{23}-c_{22}c_{23}\right) ^{2}\right.
\right.  \nonumber \\
&&+12c_{21}\left[ 3\left( c_{22}\right)
^{2}c_{14}-9c_{11}c_{14}c_{22}+c_{12}\left( c_{23}\right)
^{2}-3c_{12}c_{22}c_{24}\right.  \nonumber \\
&&\left. \left. \left. +9c_{12}c_{14}c_{21}-\left( c_{12}\right)
^{2}c_{23}-9c_{14}c_{21}c_{23}\right] \right\} ^{1/2}\right] /\left(
6c_{21}\right) ~,  \label{A3}
\end{eqnarray}

\begin{eqnarray}
c_{14} &=&\left\{ -3\left( c_{13}\right) ^{2}c_{21}-c_{12}\left[
c_{12}c_{23}+3c_{22}c_{24}-\left( c_{23}\right) ^{2}\right] \right.  \nonumber
\\
&&\left. +c_{13}\left(
c_{12}c_{22}+3c_{11}c_{23}-c_{22}c_{23}+9c_{21}c_{24}\right) \right\} /
\nonumber \\
&&/\left\{ 3\left[ 3\left( c_{11}c_{22}-c_{12}c_{21}\right)
+3c_{21}c_{23}-\left( c_{22}\right) ^{2}\right] \right\} ~.  \label{A4}
\end{eqnarray}

\section{Appendix B}

In this \textbf{Appendix B }the formulas ---implied by the \textit{second
constraint} on the $8$ coefficients $c_{n\ell }$ obtained in \textbf{%
Subsection 3.2}: see for instance eq. (\ref{32SecCona})---are reported,
which express each of these $8$ coefficients in terms of the other $7$; only
$4$ of these expressions are actually displayed below, since the other $4$
can be obtained from those displayed via the transformations (\ref{1Change}%
). These formulas provide of course $8$ \textit{constraints} on the $8$
coefficients $c_{n\ell }$. These \textit{constraints} are \textit{all}
\textit{equivalent among themselves}.
\end{subequations}
\begin{subequations}
\begin{eqnarray}
c_{11} &=&\left[ 6c_{14}c_{22}c_{23}+\left(
c_{12}c_{23}+9c_{14}c_{21}\right) \left( c_{13}-3c_{24}\right) \right.
\nonumber \\
&&\left. +\left( c_{12}c_{23}-9c_{14}c_{21}\right) R_{1}\right] /\left(
18c_{14}c_{23}\right) ~,  \label{BBc11}
\end{eqnarray}

\begin{eqnarray}
c_{12} &=&\left\{ 3\left( c_{13}\right) ^{2}c_{21}+18c_{21}\left(
c_{14}c_{23}-c_{13}c_{24}\right) -c_{13}c_{23}\left( 3c_{11}-c_{22}\right)
\right.  \nonumber \\
&&+3c_{24}\left[ c_{23}\left( 3c_{11}-c_{22}\right) +9c_{21}c_{24}\right]
\nonumber \\
&&\left. +\left[ c_{23}\left( 3c_{11}-c_{22}\right)
-3c_{13}c_{21}+9c_{21}c_{24}\right] R_{1}\right\} /\left[ 2\left(
c_{23}\right) ^{2}\right] ~,  \label{BBc12}
\end{eqnarray}

\begin{eqnarray}
&&c_{13}=\left[ c_{12}c_{23}\left( 3c_{11}-c_{22}\right) +9c_{21}\left(
3c_{11}c_{14}-c_{14}c_{22}+2c_{12}c_{24}\right) \right.  \nonumber \\
&&\left. +\left( c_{12}c_{23}-9c_{14}c_{21}\right) R_{2}\right] /\left(
6c_{12}c_{21}\right) ~,
\end{eqnarray}

\begin{eqnarray}
&&c_{14}=\left\{ -3\left( 3c_{11}-c_{22}\right) \left(
c_{13}c_{21}-3c_{21}c_{24}\right) \right.  \nonumber \\
&&+c_{23}\left[ \left( 3c_{11}-c_{22}\right) ^{2}+6c_{12}c_{21}\right]
\nonumber \\
&&\left. +\left[ 3c_{21}\left( c_{13}-3c_{24}\right) -c_{23}\left(
3c_{11}-c_{22}\right) \right] R_{2}\right\} /\left[ 54\left( c_{21}\right)
^{2}\right] ~,  \label{BBc14}
\end{eqnarray}%
with
\begin{equation}
R_{1}=\sqrt{\left( c_{13}-3c_{24}\right) ^{2}+12c_{14}c_{23}}~,~~~R_{2}=%
\sqrt{(3c_{11}-c_{22})^{2}+12c_{12}c_{21}~.}  \label{BRR}
\end{equation}

\section{Appendix C}

In this \textbf{Appendix C} we report the expressions of \textit{pairs} of
the $8$ coefficients $c_{n\ell }$ in term of the other $6$ coefficients, as
implied by the $2$ \textit{independent constraints} on these $8$
coefficients obtained in \textbf{Subsection 3.2 }(see in particular \textbf{%
Remark 3.2-2}): which are themselves \textit{sufficient} for the \textit{%
solvability} of the model (\ref{1}) as described by \textbf{Proposition 1-1}%
. These expressions---being the solutions of \textit{nonlinear} \textit{%
algebraic} equations---are generally \textit{not }unique; in some cases the
only multiplicities are those implied by the \textit{sign} ambiguities
intrinsic in the \textit{square-root} or \textit{cubic root} operations; in
other cases the multiplicities are less trivial and alternative solutions
are then displayed. We only report results for \textit{generic} values of
the coefficients $c_{n\ell }$; in particular we omit to report special
solutions with \textit{vanishing} coefficients. In all cases we were able to
obtain an \textit{explicit} solution via \textbf{Mathematica}; but in some
cases the formulas so obtained were so complicated that it made no sense to
report them here: these cases are identified below with the symbol !?!
(which might also indicate more than one such solution).

\textbf{Remark C-1}. Some reader might wonder about the usefulness of
formulas featuring this symbol "!?!". The point is that we are utilizing
here formulas featuring a lot of \textit{algebraic symbols }(such as $%
c_{nm\ell }$); in \textit{applicative} context, many of these symbols may be
replaced by \textit{numbers}, and in such cases \textbf{Mathematica}---or
other equivalent computer subroutines---is likely to yield more useful
outcomes. $\blacksquare $

There are of course $28$ ($=8\cdot 7/2$) different pairs of the $8$
coefficients $c_{n\ell }$, but only $16$ formulas are displayed below, since
the remaining $12$ formulas are then implied by those displayed via the
transformations (\ref{1Change}).

Hereafter the square-roots $R_{1}$ and $R_{2}$ are defined as at the end of
\textbf{Appendix B}. The other $2$ square-root functions featured by the
following formulas are defined as follows:
\end{subequations}
\begin{subequations}
\begin{eqnarray}
&&R_{3}=\left\{ \left[ \left( c_{11}\right) ^{2}\left(
9c_{14}c_{21}+c_{13}c_{22}-3c_{22}c_{24}\right) +\left( c_{13}c_{21}\right)
^{2}\right] ^{2}\right.  \nonumber \\
&&-4\left( c_{11}\right) ^{2}c_{13}c_{21}\left[ 9c_{13}c_{14}\left(
c_{21}\right) ^{2}+9\left( c_{11}\right) ^{2}c_{14}c_{22}+\left(
c_{13}\right) ^{2}c_{21}c_{22}\right.  \nonumber \\
&&\left. \left. -3c_{11}c_{14}\left( c_{22}\right)
^{2}-9c_{11}c_{13}c_{21}c_{24}\right] \right\} ^{1/2}~,  \label{R3}
\end{eqnarray}%
\begin{equation}
R_{4}=\sqrt{\left( c_{12}-c_{23}\right) ^{2}+4c_{13}c_{22}}~.  \label{R4}
\end{equation}

\subsection{The 4 pairs that are transformed into themselves by the
transformations (\protect\ref{1Change})}

\end{subequations}
\begin{eqnarray}
&&c_{11}=\frac{2c_{13}\left( c_{22}\right) ^{2}+\left( c_{12}-c_{23}\right)
\left( c_{12}c_{22}+3c_{13}c_{21}\right) +R_{4}\left(
c_{12}c_{22}-3c_{13}c_{21}\right) }{6c_{13}c_{22}}~,  \nonumber \\
&&c_{24}=-\left\{ 6\left( c_{13}\right) ^{3}c_{21}c_{22}+\left(
c_{13}\right) ^{2}\left\{ 3c_{21}\left( c_{23}\right) ^{2}-c_{12}\left[
2\left( c_{22}\right) ^{2}+3c_{21}c_{23}\right] \right\} \right.  \nonumber \\
&&+\left( c_{12}-c_{23}\right) \left[ 3c_{12}c_{14}\left( c_{22}\right)
^{2}+c_{13}c_{22}\left( c_{12}c_{23}-9c_{14}c_{21}\right) \right]  \nonumber \\
&&+\left. \left( c_{13}c_{23}-3c_{14}c_{22}\right) \left(
c_{12}c_{22}-3c_{13}c_{21}\right) R_{4}\right\} /\left[ 6c_{13}c_{22}\left(
c_{12}c_{22}-3c_{13}c_{21}\right) \right] ;  \nonumber \\
&&  \label{c11c24}
\end{eqnarray}%
\begin{subequations}
\begin{equation}
c_{12}=\frac{3c_{14}(3c_{11}-c_{22})}{c_{13}-3c_{24}}~,~~~c_{23}=!?!~;
\label{c12c23a}
\end{equation}%
\begin{equation}
c_{12}=!?!~,~~~c_{23}=!?!~;  \label{c12c23b}
\end{equation}%
\end{subequations}
\begin{equation}
c_{13}=!?!~,~~~c_{22}=!?!~;  \label{c13c22}
\end{equation}%
\begin{subequations}
\label{c14c21}
\begin{equation}
c_{14}=\frac{c_{12}\left( c_{13}-3c_{24}\right) }{3\left(
3c_{11}-c_{22}\right) }~,~~~c_{21}=\frac{c_{23}\left( 3c_{11}-c_{22}\right)
}{3\left( c_{13}-3c_{24}\right) }~;  \label{c14c21a}
\end{equation}%
\begin{eqnarray}
&&c_{14}=\left\{ \left( c_{12}\right) ^{2}c_{23}-c_{12}\left[ 2\left(
c_{23}\right) ^{2}+c_{22}\left( c_{13}-3c_{24}\right) \right] \right.  \nonumber
\\
&&\left. +c_{23}\left[ \left( c_{23}\right) ^{2}+3c_{22}\left(
c_{13}-c_{24}\right) \right] +FR_{4}\right\} /\left[ 6\left( c_{22}\right)
^{2}\right] ~,  \nonumber \\
&&c_{21}=\left\{ -\left( c_{12}\right) ^{4}c_{23}+c_{13}\left(
3c_{11}-c_{22}\right) \left\{ c_{23}\left[ \left( c_{23}\right)
^{2}+c_{22}\left( c_{13}-3c_{24}\right) \right] -FR_{4}\right\} \right.
\nonumber \\
&&+\left( c_{12}\right) ^{3}\left[ 3\left( c_{23}\right) ^{2}+c_{22}\left(
c_{13}-3c_{24}\right) \right] +c_{12}\left\{ -3\left( c_{13}\right)
^{2}c_{22}\left( c_{11}-c_{22}\right) +\left( c_{23}\right) ^{4}\right.
\nonumber \\
&&\left. -c_{23}FR_{4}-3c_{22}\left( c_{23}\right) ^{2}c_{24}+c_{13}\left[
9c_{22}c_{24}\left( c_{11}-c_{22}\right) -\left( c_{23}\right) ^{2}\left(
6c_{11}-5c_{22}\right) \right] \right\}  \nonumber \\
&&\left. +\left( c_{12}\right) ^{2}\left\{ -3\left( c_{23}\right)
^{3}+FR_{4}+c_{23}\left[ 3c_{11}c_{13}-c_{22}\left( 5c_{13}-6c_{24}\right) %
\right] \right\} \right\} /\left[ 6F\left( c_{13}\right) ^{2}\right] ~,
\nonumber \\
&&F=c_{22}\left( c_{13}-3c_{24}\right) -c_{23}\left( c_{12}-c_{23}\right) ~.
\label{c14c21b}
\end{eqnarray}

\subsection{The 12 pairs that are transformed into 12 \textit{other} pairs
by the transformations (\protect\ref{1Change})}

\end{subequations}
\begin{subequations}
\label{C11C12}
\begin{equation}
c_{11}=\frac{c_{22}c_{23}+3c_{21}\left( c_{13}-3c_{24}\right) }{3c_{23}}%
~,~~~c_{12}=\frac{9c_{14}c_{21}}{c_{23}}~;  \label{c11c12a}
\end{equation}%
\begin{eqnarray}
&&c_{11}=\left\{ \left( c_{13}\right) ^{2}c_{23}\left( c_{13}-6c_{24}\right)
+27\left( c_{14}\right) ^{2}c_{21}\left( c_{13}-3c_{24}\right) \right.
\nonumber \\
&&+9c_{23}\left[ c_{13}\left( c_{24}\right) ^{2}+c_{14}c_{23}\left(
c_{13}-c_{24}\right) \right]  \nonumber \\
&&\left. -\left\{ c_{13}c_{23}\left( c_{13}-3c_{24}\right) +3c_{14}\left[
\left( c_{23}\right) ^{2}-9c_{14}c_{21}\right] \right\} R_{1}\right\} /\left[
\left( 54\left( c_{14}\right) ^{2}c_{23}\right) \right] ~,  \nonumber \\
&&c_{12}=\left\{ c_{13}\left(
c_{13}c_{23}+3c_{14}c_{22}-3c_{23}c_{24}\right) +3c_{14}\left[ 2\left(
c_{23}\right) ^{2}-3c_{22}c_{24}\right] \right.  \nonumber \\
&&\left. +\left( 3c_{14}c_{22}-c_{13}c_{23}\right) R_{1}\right\} /\left(
6c_{14}c_{23}\right) ~;  \label{c11c12b}
\end{eqnarray}%
\end{subequations}
\begin{eqnarray}
&&c_{11}=\left\{ 9\left( c_{14}\right) ^{2}c_{22}\left[ 3c_{12}c_{21}+\left(
c_{22}\right) ^{2}-6c_{21}c_{23}\right] \right.  \nonumber \\
&&+3\left( c_{12}\right) ^{2}c_{23}\left( c_{23}c_{24}-c_{14}c_{22}\right)
+3c_{14}c_{22}\left( c_{23}\right) ^{2}\left( c_{12}+c_{23}\right)  \nonumber \\
&&-9c_{14}\left( c_{22}\right) ^{2}c_{24}\left( c_{12}+c_{23}\right)
-3c_{12}c_{23}c_{24}\left[ \left( c_{23}\right) ^{2}-3c_{22}c_{24}\right]
\nonumber \\
&&+C_{13}\left[ 3c_{14}c_{22}\left( c_{12}c_{22}+9c_{21}c_{24}\right)
-c_{12}c_{23}\left( c_{12}c_{23}+6c_{22}c_{24}\right) \right.  \nonumber \\
&&\left. \left. +9c_{14}c_{21}\left( c_{23}\right) ^{2}\right] +\left(
C_{13}\right) ^{2}c_{22}\left( c_{12}c_{23}-9c_{14}c_{21}\right) \right\} /
\nonumber \\
&&/\left\{ 9c_{14}\left[ 3c_{22}\left( c_{14}c_{22}-c_{23}c_{24}\right)
+\left( c_{23}\right) ^{2}\left( c_{23}-c_{12}\right) \right] \right\} ~,
\nonumber \\
&&c_{13}\equiv C_{13}=2c_{24}+\left[
c_{12}c_{23}+2^{1/3}A_{3}/A_{4}-2^{-1/3}A_{4}\right] /\left( 3c_{22}\right)
~,  \nonumber \\
&&A_{2}=27c_{14}c_{22}\left[ -9c_{14}\left( c_{22}\right) ^{3}+2\left(
c_{12}\right) ^{2}c_{22}c_{23}-3c_{12}c_{22}\left( c_{23}\right) ^{2}\right.
\nonumber \\
&&\left. +3c_{22}\left( c_{23}\right) ^{3}+3\left( c_{22}\right)
^{2}c_{24}\left( c_{12}+3c_{23}\right) \right] -9c_{12}c_{22}\left(
c_{23}\right) ^{2}c_{24}\left( c_{12}+3_{23}\right)  \nonumber \\
&&+27\left( c_{22}\right) ^{2}c_{23}\left( c_{24}\right) ^{2}\left(
c_{12}-6c_{23}\right) +54\left( c_{22}c_{24}\right) ^{3}-2\left(
c_{12}c_{23}\right) ^{3}~,  \nonumber \\
&&A_{3}=9c_{14}\left( c_{22}\right) ^{2}\left( 3c_{23}-c_{12}\right)
-c_{12}c_{23}\left( c_{12}c_{23}+3c_{22}c_{24}\right)  \nonumber \\
&&-9c_{22}c_{24}\left[ \left( c_{23}\right) ^{2}+c_{22}c_{24}\right] ~,
\nonumber \\
&&A_{4}=\left( A_{2}+\sqrt{\left( A_{2}\right) ^{2}+4\left( A_{3}\right) ^{3}%
}\right) ^{1/3}~;  \label{c11c13}
\end{eqnarray}

\begin{subequations}
\begin{equation}
c_{11}=\frac{3c_{13}c_{21}+c_{22}c_{23}-9c_{21}c_{24}}{3c_{23}}~,~~~c_{14}=%
\frac{c_{12}c_{23}}{9c_{21}}~;  \label{c11c14a}
\end{equation}%
\begin{eqnarray}
&&c_{11}=\left\{ -\left( c_{12}\right) ^{3}c_{22}c_{23}+c_{13}\left[ 2\left(
c_{22}\right) ^{2}-3c_{21}c_{23}\right] \left[ c_{13}c_{22}+\left(
c_{23}\right) ^{2}-3c_{22}c_{24}\right] \right.  \nonumber \\
&&+\left( c_{12}\right) ^{2}\left\{ c_{13}\left[ \left( c_{22}\right)
^{2}-3c_{21}c_{23}\right] +c_{22}\left[ 2\left( c_{23}\right)
^{2}-3c_{22}c_{24}\right] \right\}  \nonumber \\
&&+c_{12}\left\{ 3\left( c_{13}\right) ^{2}c_{21}c_{22}-c_{22}c_{23}\left[
\left( c_{23}\right) ^{2}-3c_{22}c_{24}\right] \right.  \nonumber \\
&&\left. -3c_{13}\left[ c_{23}\left[ \left( c_{22}\right) ^{2}-2c_{21}c_{23}%
\right] +3c_{21}c_{22}c_{24}\right] \right\}  \nonumber \\
&&\left. +\left( 3c_{13}c_{21}-c_{12}c_{22}\right) \left[
c_{13}c_{22}-c_{12}c_{23}+\left( c_{23}\right) ^{2}-3c_{22}c_{24}\right]
R_{4}\right\} /  \nonumber \\
&&\left\{ 6c_{13}c_{22}\left[ c_{13}c_{22}-c_{12}c_{23}+\left( c_{23}\right)
^{2}-3c_{22}c_{24}\right] \right\} ~,  \nonumber \\
&&c_{14}=\left\{ c_{12}c_{23}\left( c_{12}-2c_{23}\right) +3c_{22}\left[
c_{12}c_{24}+c_{23}\left( c_{13}-c_{24}\right) \right] -c_{12}c_{13}c_{22}%
\right.  \nonumber \\
&&\left. +\left( c_{23}\right) ^{3}+\left[ c_{13}c_{22}-c_{12}c_{23}+\left(
c_{23}\right) ^{2}-3c_{22}c_{24}\right] R_{4}\right\} /\left[ 6\left(
c_{22}\right) ^{2}\right] ~;  \label{c11c14b}
\end{eqnarray}

\end{subequations}
\begin{equation}
c_{11}=\frac{c_{12}c_{13}+3c_{14}c_{22}-3c_{12}c_{24}}{9c_{14}}~,~~~c_{21}=%
\frac{c_{12}c_{23}}{9c_{14}}~;  \label{c11c21}
\end{equation}%
\begin{eqnarray}
&&c_{11}=\left[ 54\left( c_{14}\right) ^{2}c_{23}\right] ^{-1}\left\{ \left(
c_{13}\right) ^{3}c_{23}-6\left( c_{13}\right) ^{2}c_{23}c_{24}\right.
\nonumber \\
&&+9c_{13}\left[ 3\left( c_{14}\right) ^{2}c_{21}+c_{14}\left( c_{23}\right)
^{2}+c_{23}\left( c_{24}\right) ^{2}\right]  \nonumber \\
&&-9c_{14}\left[ \left( c_{23}\right) ^{2}c_{24}+9c_{14}c_{21}c_{24}\right]
\nonumber \\
&&\left. +\left[ -\left( c_{13}\right)
^{2}c_{23}+3c_{13}c_{23}c_{24}-3c_{14}\left( c_{23}\right) ^{2}+27\left(
c_{14}\right) ^{2}c_{21}\right] R_{1}\right\} ~,  \nonumber \\
&&c_{22}=\left[ 18\left( c_{14}\right) ^{2}\right] ^{-1}\left\{ \left(
c_{13}\right) ^{3}+3c_{14}\left[
-c_{12}c_{13}+3c_{13}c_{23}+3c_{12}c_{24}-3c_{23}c_{24}\right] \right.
\nonumber \\
&&\left. +3c_{13}c_{24}\left( 3c_{24}-2c_{13}\right) +\left[ c_{13}\left(
c_{13}-3c_{24}\right) -3c_{14}\left( c_{12}-c_{23}\right) \right]
R_{1}\right\} ~;  \label{c11c22}
\end{eqnarray}%
\begin{subequations}
\label{c11c23}
\begin{equation}
c_{11}=\frac{c_{12}c_{13}+3\left( c_{14}c_{22}-c_{12}c_{24}\right) }{9c_{14}}%
~,~~~c_{23}=\frac{9c_{14}c_{21}}{c_{12}}~;  \label{c11c23a}
\end{equation}%
\begin{eqnarray}
&&c_{11}=\left\{ 27\left( c_{13}c_{14}\right) ^{2}c_{22}\left[ c_{13}\left(
3c_{13}c_{21}-c_{12}c_{22}-9c_{21}c_{24}\right) \right. \right.  \nonumber \\
&&\left. +3c_{12}\left( c_{22}c_{24}-3c_{14}c_{21}\right) \right] +81\left(
c_{14}\right) ^{3}\left( c_{22}\right) ^{2}\cdot  \nonumber \\
&&\left. \left[ \left( c_{12}\right) ^{2}c_{13}+3c_{14}\left(
3c_{13}c_{21}-c_{12}c_{22}\right) \right] \right\} ^{-1}\cdot  \nonumber \\
&&\cdot \left\{ 9\left( c_{12}c_{13}c_{14}\right) ^{2}\left[ c_{12}\left(
c_{22}\right) ^{2}+27c_{21}c_{22}c_{24}-27c_{14}\left( c_{21}\right)
^{2}-9c_{13}c_{21}c_{22}\right] \right.  \nonumber \\
&&+27\left( c_{13}c_{14}\right) ^{2}c_{21}\left\{ 3c_{12}\left[ \left(
c_{13}\right) ^{2}c_{21}+c_{14}\left( c_{22}\right) ^{2}-3c_{13}c_{21}c_{24}%
\right] \right.  \nonumber \\
&&\left. +c_{22}\left[ \left( c_{13}\right) ^{2}c_{22}+3\left(
9c_{14}c_{21}-c_{13}c_{22}\right) c_{24}\right] \right\}  \nonumber \\
&&+81\left( c_{14}\right) ^{3}\left( c_{22}\right) ^{2}\left[
c_{14}c_{22}\left( 3c_{13}c_{21}-c_{12}c_{22}\right) +\left( c_{12}\right)
^{2}\left( c_{22}c_{24}-3c_{14}c_{21}\right) \right]  \nonumber \\
&&+9\left( c_{12}c_{13}c_{22}\right) ^{2}c_{14}c_{24}\left(
2c_{13}-3c_{24}\right)  \nonumber \\
&&+9c_{12}\left( c_{13}\right) ^{3}c_{14}c_{22}\left[ 2\left( c_{13}\right)
^{2}c_{21}-c_{14}\left( c_{22}\right)
^{2}-12c_{13}c_{21}c_{24}+18c_{21}\left( c_{24}\right) ^{2}\right]  \nonumber \\
&&+27\left( c_{12}\right) ^{3}c_{13}\left( c_{14}\right) ^{2}c_{22}\left(
3c_{14}c_{21}-c_{22}c_{24}\right)  \nonumber \\
&&+27c_{12}c_{13}\left( c_{14}\right) ^{2}c_{22}\left[ 27\left(
c_{14}c_{21}\right) ^{2}-c_{22}c_{24}\left(
18c_{14}c_{21}-c_{13}c_{22}\right) \right]  \nonumber \\
&&-3\left( c_{13}\right) ^{3}c_{14}c_{22}\left[ \left( c_{12}\right)
^{2}c_{13}c_{22}+162\left( c_{14}c_{21}\right) ^{2}\right]  \nonumber \\
&&+C_{23}\left\{ 3\left( c_{12}c_{13}\right) ^{2}c_{14}\left[ 6\left(
c_{13}\right)
^{2}c_{21}+18c_{12}c_{14}c_{21}+2c_{12}c_{13}c_{22}-3c_{14}\left(
c_{22}\right) ^{2}\right. \right.  \nonumber \\
&&+\left. 6\left( 3c_{13}c_{21}-c_{12}c_{22}\right) c_{24}\right] +\left(
c_{12}\right) ^{2}\left( c_{13}\right) ^{3}c_{22}\left[ 6c_{13}c_{24}-9%
\left( c_{24}\right) ^{2}-\left( c_{13}\right) ^{2}\right]  \nonumber \\
&&+27\left( c_{14}c_{22}\right) ^{2}\left[ \left( c_{12}\right)
^{3}c_{14}+9c_{12}\left( c_{14}\right) ^{2}c_{21}+\left( c_{12}\right)
^{2}c_{13}c_{14}\right]  \nonumber \\
&&+27\left( c_{13}c_{14}\right) ^{2}c_{21}\left(
18c_{12}c_{14}c_{21}+5c_{12}c_{13}c_{22}+9c_{13}c_{21}c_{24}\right)  \nonumber
\\
&&-243c_{13}\left( c_{14}\right) ^{3}c_{21}c_{22}\left[ \left( c_{12}\right)
^{2}+3c_{14}c_{21}\right]  \nonumber \\
&&\left. -9c_{12}c_{13}\left( c_{14}\right) ^{2}c_{22}\left[ \left(
c_{12}\right) ^{3}+27c_{13}c_{21}c_{24}\right] \right\}  \nonumber \\
&&\left( C_{23}\right) ^{2}\left\{ \left( c_{12}c_{13}\right) ^{3}\left(
c_{13}-3c_{24}\right) +9\left( c_{12}c_{13}\right) ^{2}c_{14}\left[
c_{22}\left( 2c_{24}-c_{13}\right) -12c_{14}c_{21}\right] \right.  \nonumber \\
&&+81c_{13}\left( c_{14}\right) ^{3}c_{21}\left[ 2c_{12}c_{22}-3c_{13}c_{21}%
\right] -54c_{12}\left( c_{13}\right) ^{3}c_{14}c_{21}c_{24}  \nonumber \\
&&+\left. 3\left( c_{12}\right) ^{2}c_{14}\left[ 6c_{12}c_{13}c_{14}c_{22}-9%
\left( c_{14}c_{22}\right) ^{2}-\left( c_{12}c_{13}\right) ^{2}\right]
\right\}  \nonumber \\
&&+3\left( C_{23}\right) ^{3}\left[ \left( c_{12}c_{13}\right) ^{2}\left(
2c_{12}c_{14}+c_{13}c_{24}\right) +3c_{12}c_{13}\left( c_{14}\right)
^{2}\left( 6c_{13}c_{21}-c_{12}c_{22}\right) \right]  \nonumber \\
&&\left. -3\left( C_{23}\right) ^{4}\left( c_{12}c_{13}\right)
^{2}c_{14}\right\} ~,  \nonumber \\
&&c_{23}\equiv C_{23}=2c_{12}/3+\left[ 3c_{13}c_{24}-2^{1/3}\left(
E_{2}/E_{3}\right) +2^{-1/3}E_{3}\right] /\left( 9c_{14}\right) ~,  \nonumber
\end{eqnarray}%
\begin{eqnarray}
&&E_{1}=27\left\{ \left( c_{12}c_{14}\right) ^{2}\left[ 6\left(
c_{13}\right) ^{2}-2c_{12}c_{14}-3c_{13}c_{24}\right] +9\left(
c_{13}c_{14}\right) ^{2}c_{22}\left( 3c_{24}-c_{13}\right) \right.  \nonumber \\
&&+3c_{12}c_{14}c_{24}\left[ \left( c_{13}\right) ^{3}-9\left( c_{14}\right)
^{2}c_{22}\right] +27\left( c_{14}\right) ^{3}c_{22}\left(
3c_{14}c_{22}-c_{12}c_{13}\right)  \nonumber \\
&&\left. +c_{13}\left( c_{24}\right) ^{2}\left[ 3c_{12}c_{13}c_{14}-45\left(
c_{14}\right) ^{2}c_{22}+2\left( c_{13}\right) ^{2}c_{24}\right] \right\} ~,
\nonumber \\
&&E_{2}=9\left\{ \left( c_{14}\right) ^{2}\left[ 9c_{22}\left(
c_{13}-c_{24}\right) -\left( c_{12}\right) ^{2}\right] -\left(
c_{13}c_{24}\right) ^{2}-c_{12}c_{13}c_{14}\left( c_{13}+c_{14}\right)
\right\} ~,  \nonumber \\
&&E_{3}=\left( E_{1}+\sqrt{\left( E_{1}\right) ^{2}+4\left( E_{2}\right) ^{3}%
}\right) ^{1/3}~.
\end{eqnarray}

\end{subequations}
\begin{subequations}
\label{c12c13}
\begin{equation}
c_{12}=\frac{9c_{14}c_{21}}{c_{23}}~,~~~c_{13}=\frac{%
3c_{11}c_{23}-c_{22}c_{23}+9c_{21}c_{24}}{3c_{21}}~;  \label{c12c13a}
\end{equation}%
\begin{equation}
c_{12}=!?!~,~~~c_{13}=!?!  \label{c12c13b}
\end{equation}%
\end{subequations}
\begin{eqnarray}
&&c_{12}\equiv C_{12}=\left( 2/3\right) c_{21}c_{23}+\left(
3c_{11}c_{22}+2^{1/3}B_{3}/B_{4}-2^{-1/3}B_{4}\right) /\left( 9c_{21}\right)
~,  \nonumber \\
&&c_{14}=\left\{ 3\left( c_{13}\right) ^{2}c_{21}\left(
3c_{11}-2c_{22}\right) -9c_{11}c_{13}c_{23}\left( c_{11}-c_{22}\right)
\right.  \nonumber \\
&&-2c_{13}\left( c_{22}\right) ^{2}c_{23}-9c_{13}c_{21}c_{24}\left(
3c_{11}-2c_{22}\right)  \nonumber \\
&&+C_{12}\left[ \left( c_{22}\right) ^{2}\left( c_{13}-3c_{24}\right)
-\left( c_{23}\right) ^{2}\left( 3c_{11}-2c_{22}\right) -9c_{21}c_{23}c_{24}%
\right]  \nonumber \\
&&\left. +\left( C_{12}\right) ^{2}\left[ -3c_{21}\left(
c_{13}-3c_{24}\right) +c_{23}\left( 3c_{11}-c_{22}\right) \right] \right\} /
\nonumber \\
&&/\left\{ 3\left[ \left( c_{22}\right) ^{2}\left( 3c_{11}-c_{22}\right)
+3c_{21}\left( c_{22}c_{23}-3c_{13}c_{21}\right) \right] \right\} ~,  \nonumber
\\
&&B_{2}=27c_{13}\left( c_{21}\right) ^{2}\left[ -81c_{13}\left(
c_{21}\right) ^{2}+54\left( c_{11}\right) ^{2}c_{22}-27c_{11}\left(
c_{22}\right) ^{2}\right.  \nonumber \\
&&\left. +9\left( c_{22}\right) ^{3}+27c_{21}c_{23}\left(
c_{11}+c_{22}\right) \right] -2\left[ \left( c_{11}c_{22}\right) ^{3}-\left(
c_{21}c_{23}\right) ^{3}\right]  \nonumber \\
&&-3c_{11}c_{21}c_{22}c_{23}\left[ c_{11}c_{22}+\left( c_{22}\right)
^{2}-c_{21}c_{23}\right] -6\left( c_{21}c_{22}c_{23}\right) ^{2}~,  \nonumber \\
&&B_{3}=-81c_{13}\left( c_{21}\right) ^{2}\left( c_{11}-c_{22}\right)
-9c_{21}c_{22}c_{23}\left( c_{11}+c_{22}\right)  \nonumber \\
&&-9\left[ \left( c_{11}c_{22}\right) ^{2}+\left( c_{21}c_{23}\right) ^{2}%
\right] ~,  \nonumber \\
&&B_{4}=\left( B_{2}+\sqrt{\left( B_{2}\right) ^{2}+4\left( B_{3}\right) ^{3}%
}\right) ^{1/3}~.  \label{c12c14aa}
\end{eqnarray}%
\begin{subequations}
\begin{equation}
c_{12}=\frac{9c_{11}c_{14}-3c_{14}c_{22}}{c_{13}-3c_{24}}~,~~~c_{21}=\frac{%
3c_{11}c_{23}-c_{22}c_{23}}{3\left( c_{13}-3c_{24}\right) }~;
\label{c12c21a}
\end{equation}%
\begin{eqnarray}
&&c_{12}=\left[ \left( c_{13}-3c_{24}\right) \left(
c_{13}c_{23}+3c_{14}c_{22}\right) \right.  \nonumber \\
&&\left. +6c_{14}\left( c_{23}\right) ^{2}+\left(
3c_{14}c_{22}-c_{13}c_{23}\right) R_{1}\right] /\left( 6c_{14}c_{23}\right)
~,  \nonumber \\
&&c_{21}=\left\{ \left( c_{13}\right) ^{3}\left( c_{13}-9c_{24}\right)
+9\left( c_{14}\right) ^{2}\left[ 2\left( c_{23}\right) ^{2}+9c_{11}c_{24}%
\right] \right.  \nonumber \\
&&+3c_{14}c_{23}\left[ 9\left( c_{24}\right) ^{2}+4\left( c_{13}\right) ^{2}%
\right] +3\left( c_{13}\right) ^{2}\left[ 4c_{14}c_{23}+9\left(
c_{24}\right) ^{2}\right]  \nonumber \\
&&-9c_{13}\left[ 3c_{11}\left( c_{14}\right) ^{2}+3\left( c_{24}\right)
^{3}+5c_{14}c_{23}c_{24}\right]  \nonumber \\
&&+\left[ -\left( c_{13}\right) ^{3}+6\left( c_{13}\right)
^{2}c_{24}-9c_{13}\left( c_{24}\right) ^{2}-6c_{13}c_{14}c_{23}\right.
\nonumber \\
&&\left. \left. +9c_{14}\left( 3c_{11}c_{14}+c_{23}c_{24}\right) \right]
R_{1}\right\} /\left[ 162\left( c_{14}\right) ^{3}\right] ~;  \label{c12c21b}
\end{eqnarray}%
\end{subequations}
\begin{equation}
c_{12}=\frac{9c_{14}c_{21}}{c_{23}}~,~~~c_{22}=\frac{%
3(c_{11}c_{23}-c_{13}c_{21}+3c_{21}c_{24})}{c_{23}}~;  \label{c12c22}
\end{equation}%
\begin{subequations}
\begin{equation}
c_{13}=\frac{3c_{11}c_{23}-c_{22}c_{23}}{3c_{21}}~,~~~c_{14}=\frac{%
c_{12}c_{23}}{9c_{21}}~;  \label{c13c14a}
\end{equation}%
\begin{eqnarray}
&&c_{13}=\left\{ c_{22}\left[ 9\left( c_{11}\right)
^{2}+9c_{12}c_{21}+\left( c_{22}\right) ^{2}-3c_{21}c_{23}\right] \right.
\nonumber \\
&&-3c_{11}\left[ 3c_{12}c_{21}+2\left( c_{22}\right) ^{2}-c_{21}c_{23}\right]
\nonumber \\
&&\left. +\left[ 3c_{21}\left( c_{12}-c_{23}\right) +c_{22}\left(
3c_{11}+c_{22}\right) \right] R_{2}\right\} /\left[ 18\left( c_{21}\right)
^{2}\right] ~,  \nonumber \\
&&c_{14}=\left\{ 27\left( c_{11}\right) ^{2}\left[
c_{12}c_{21}-c_{11}c_{22}+\left( c_{22}\right) ^{2}\right]
+9c_{12}c_{21}\left( 2c_{12}c_{21}-5c_{11}c_{22}\right) \right.  \nonumber \\
&&+\left( c_{22}\right) ^{2}\left[ 12c_{12}c_{21}-9c_{11}c_{22}+\left(
c_{22}\right) ^{2}\right] +\left( c_{21}\right) ^{2}c_{24}\left(
81c_{11}-27c_{22}\right)  \nonumber \\
&&-\left[ 9c_{21}\left( c_{11}c_{12}+3c_{21}c_{24}\right) -c_{22}\left[
6c_{12}c_{21}+\left( c_{22}\right) ^{2}\right] \right.  \nonumber \\
&&\left. \left. +3c_{11}c_{22}\left( 2c_{22}-3c_{11}\right) \right]
R_{2}\right\} /\left[ 162\left( c_{21}\right) ^{3}\right] ~;  \label{c13c14b}
\end{eqnarray}%
\end{subequations}
\begin{subequations}
\begin{equation}
c_{13}=\frac{3(3c_{11}c_{14}+c_{12}c_{24}-c_{14}c_{22})}{c_{12}}~,~~~c_{21}=%
\frac{c_{12}c_{23}}{9c_{14}}~;  \label{c13c21a}
\end{equation}%
\begin{equation}
c_{13}=!?!~,~~~c_{21}=!?!~;  \label{c13c21b}
\end{equation}

\end{subequations}


\bigskip

\bigskip









\begin{thebibliography}{}

\bibitem{GRN2016} G. R. Nicklason, "Homogeneous-Like Generalized Cubic
Systems", Internat. J. Diff. Eqs., Vol. 2016, Article ID 7640340 (15 pages);
https://dx.doi.org/10.1155/2016/7640340.

\bibitem{CP2019} F. Calogero and F. Payandeh, \textquotedblleft Polynomials
with multiple zeros and solvable dynamical systems including models in the
plane with polynomial interactions\textquotedblright , J. Math. Phys.
\textbf{60}, 082701 (2019) (23 pages); doi.org/101063.1.5082249;
arXiv:1904.00496v1 [math-ph] 31 Mar 2019.

\bibitem{CCL2020} F. Calogero, R. Conte and F. Leyvraz, \textquotedblleft
New algebraically solvable systems of two autonomous first-order ordinary
differential equations with purely quadratic right-hand
sides\textquotedblright , J. Math. Phys. \textbf{61}, 102704 (2020); 
https//doi.org/10.1063/5.0011257; http//arxive.org/abs/2009.11200.

\bibitem{RG1960} R. Garnier, "Sur des syst\`{e}mes diff\'{e}rentielles du
second ordre dont l'int\'{e}gral general est uniform", Ann. \'{E}cole Norm.
\textbf{77}(2) 123-144 (1960).

\bibitem{FC2008} F. Calogero, \textit{Isochronous Systems}, Oxford
University Press, Oxford, U. K. (264 pages), hardback (2008), updated
paperback (2012).
\end{thebibliography}
\end{document}